%% file: main.tex
\documentclass{article}
\usepackage[utf8]{inputenc}
\usepackage{zibtitlepage}
\usepackage[ruled,vlined]{algorithm2e}
\usepackage{cite}
\usepackage{xspace}
\usepackage{amssymb,amsthm,amsmath}
\usepackage{amsfonts}
\usepackage[table]{xcolor}    
\usepackage{graphicx}
\graphicspath{ {./sections/Pictures/} }
\usepackage{subfig}
\usepackage{parskip}
\usepackage{tablefootnote}
\usepackage{threeparttable}
\usepackage{hyperref}
\hypersetup{colorlinks,linkcolor={blue},citecolor={blue},urlcolor={red}}

\usepackage{accents}
\usepackage{interval}
\usepackage{booktabs}
\usepackage[bottom]{footmisc}  
\newtheorem{definition}{Definition}
\input{commands}
\begin{document}

\ZTPAuthor{\ZTPHasOrcid{Lovis Anderson}{0000-0002-4316-1862},
           \ZTPHasOrcid{Mark Turner}{0000-0001-7270-1496},
           \ZTPHasOrcid{Thorsten Koch}{0000-0002-1967-0077}}

\ZTPTitle{Generative deep learning for decision making in gas networks}
\ZTPNumber{20-38}
\ZTPMonth{December}
\ZTPYear{2020}

\title{Generative deep learning for decision making in gas networks}
\author{Lovis Anderson, Mark Turner, Thorsten Koch}
\zibtitlepage

\maketitle

\begin{abstract}
\input{sections/Abstract}
\end{abstract}

\section{Introduction}
\input{sections/Introduction}

\section{Background and Related Work}
\input{sections/RelatedWork}

\section{The Solution Framework}
\input{sections/NeuralNetworkDesign}

\section{The Gas Transport Model} \label{sec:The Gas Transport Model}
\input{sections/GasTransportModel}

\section{Computational Experiments} \label{sec:Computational Experiments}
\input{sections/ComputationalExperiments}

\section{Computational Results} \label{sec:Computation Results}
\input{sections/Results}

\section{Conclusion}
\input{sections/Conclusion}

\section*{Acknowledgements}
\input{sections/Funding}

\bibliographystyle{abbrv}
\bibliography{main.bib}
\appendix
\input{sections/appendix/appendix}
\end{document}

%% file: commands.tex
\newcommand{\ubar}[1]{\underaccent{\bar}{#1}}

\newcommand{\ub}[1]{\ensuremath{\bar{#1}}\xspace}
\newcommand{\lb}[1]{\ensuremath{\ubar{#1}}\xspace}

\newcommand{\graph}{\ensuremath{G}\xspace}
\newcommand{\generalGraph}{\ensuremath{\graph=(\setVertices,\setArcs)}\xspace}

\newcommand{\setVertices}{\ensuremath{\mathcal{V}}\xspace}
\newcommand{\setVerticesAll}{\ensuremath{\setVertices = \setBoundaryNodes \cup \setInnerNodes}\xspace}
\newcommand{\setEntries}{\ensuremath{\setVertices^{+}}\xspace}
\newcommand{\setExits}{\ensuremath{\setVertices^{-}}\xspace}
\newcommand{\setBoundaryNodes}{\ensuremath{\setVertices^\text{b}}\xspace}
 
\newcommand{\setInnerNodes}{\ensuremath{\setVertices^{0}}\xspace}
\newcommand{\setFenceNodes}{\ensuremath{\setVertices^{\text{fn}}}\xspace}
\newcommand{\setFenceNodesNaviStation}[1]{\ensuremath{\setFenceNodes_{i}}\xspace}

\newcommand{\setArcs}{\ensuremath{\mathcal{A}}\xspace}
\newcommand{\setArcsAllStationModel}{\ensuremath{\setArcs = \setPipes \cup \setValves \cup \setResistors \cup \setRegulators \cup \setCompressorStations}\xspace}
\newcommand{\setPipes}{\ensuremath{\setArcs^{\text{pi}}}\xspace}
\newcommand{\setRegulators}{\ensuremath{\setArcs^{\text{rg}}}\xspace}
\newcommand{\setValves}{\ensuremath{\setArcs^{\text{va}}}\xspace}

\newcommand{\setResistors}{\ensuremath{\setArcs^{\text{rs}}}\xspace}
\newcommand{\setCompressorStations}{\ensuremath{\setArcs^{\text{cs}}}\xspace}

\newcommand{\setCompressorConfigurations}[1]{\ensuremath{\mathcal{C}_{#1}}}

\newcommand{\setConfigurationFacets}[1]{\ensuremath{\mathcal{H}_{#1}}}

\newcommand{\setNSModes}{\ensuremath{\mathcal{O}}\xspace}
\newcommand{\setNSFenceGroups}{\ensuremath{\mathcal{G}}\xspace}

\newcommand{\setTimesteps}{\ensuremath{\mathcal{T}_0}\xspace}
\newcommand{\setTimestepsNoZero}{\ensuremath{\mathcal{T}}\xspace}
\newcommand{\nTimesteps}{\ensuremath{k}\xspace}
\newcommand{\setTimestepsAll}{\ensuremath{\setTimesteps := \{0, \dots, \nTimesteps\}}\xspace}
\newcommand{\setTimestepsNoZeroAll}{\ensuremath{\setTimestepsNoZero := \setTimesteps\setminus\{0\}}\xspace}
\newcommand{\generalGranularity}{\ensuremath{\tau}\xspace}
\newcommand{\granularity}[1]{\ensuremath{\generalGranularity(#1)}\xspace}

\newcommand{\generalInflow}{\ensuremath{d}}

\newcommand{\varInflowValue}[2]{\ensuremath{\generalInflow_{#1,#2}}\xspace}

\newcommand{\demandPressure}[2]{\ensuremath{\hat{\generalPressure}_{#1,#2}}\xspace}
\newcommand{\demandInflow}[2]{\ensuremath{\hat{\generalInflow}_{#1,#2}}\xspace}

\newcommand{\generalPressure}{\ensuremath{p}\xspace}
\newcommand{\varPressure}[2]{\ensuremath{\generalPressure_{#1,#2}}\xspace}

\newcommand{\paramPressureUB}[2]{\ensuremath{\ub{\generalPressure}_{#1,#2}}\xspace}
\newcommand{\paramPressureLB}[2]{\ensuremath{\lb{\generalPressure}_{#1,#2}}\xspace}

\newcommand{\generalFlow}{\ensuremath{q}\xspace}
\newcommand{\varArcFlow}[2]{\ensuremath{\generalFlow_{#1,#2}}\xspace}

\newcommand{\paramArcFlowUB}[2]{\ensuremath{\ub{\generalFlow}_{#1,#2}}\xspace}
\newcommand{\paramArcFlowLB}[2]{\ensuremath{\lb{\generalFlow}_{#1,#2}}\xspace} 

\newcommand{\varPipeFlow}[3]{\ensuremath{\generalFlow_{#1,#2,#3}}\xspace}

\newcommand{\setFlowDirectionConditions}{\ensuremath{\mathcal{W}}\xspace}
\newcommand{\setFlowDirectionConditionsNaviStation}[1]{\ensuremath{\setFlowDirectionConditions_{i}}\xspace}

\newcommand{\generalMode}{\ensuremath{m}}
\newcommand{\varModeOpen}[2]{\ensuremath{\generalMode_{#1,#2}^\text{op}}}
\newcommand{\varModeClosed}[2]{\ensuremath{\generalMode_{#1,#2}^\text{cl}}}
\newcommand{\varModeBypass}[2]{\ensuremath{\generalMode_{#1,#2}^\text{by}}}

\newcommand{\varModeConfiguration}[3]{\ensuremath{\generalMode_{#1,#2,#3}^\text{cf}}}
\newcommand{\varNSmode}[2]{\ensuremath{\generalMode_{#1,#2}^\text{om}}}

\newcommand{\varConfigPressureL}[3]{\ensuremath{\generalPressure_{#1,#2,#3}^\text{u-cf}}\xspace}
\newcommand{\varConfigPressureR}[3]{\ensuremath{\generalPressure_{#1,#2,#3}^\text{v-cf}}\xspace}
\newcommand{\varConfigFlow}[3]{\ensuremath{\generalFlow_{#1,#2,#3}^\text{cf}}\xspace}

\newcommand{\generalSlack}{\ensuremath{\sigma}}
\newcommand{\varSlackPressurePos}[2]{\ensuremath{\generalSlack_{#1,#2}^{\generalPressure+}}}
\newcommand{\varSlackPressureNeg}[2]{\ensuremath{\generalSlack_{#1,#2}^{\generalPressure-}}}
\newcommand{\varSlackFlowPos}[2]{\ensuremath{\generalSlack_{#1,#2}^{\generalInflow+}}}
\newcommand{\varSlackFlowNeg}[2]{\ensuremath{\generalSlack_{#1,#2}^{\generalInflow-}}}

\newcommand{\generalVelocity}{\ensuremath{v}\xspace}
\newcommand{\paramAbsoluteVelocity}[2]{\ensuremath{|\generalVelocity_{#1,#2}|}\xspace} 
\newcommand{\gasConst}{\ensuremath{R}\xspace}
\newcommand{\specificGasConstant}{\ensuremath{\gasConst_s}\xspace}
\newcommand{\generalGasTemperature}{\ensuremath{T}\xspace}

\newcommand{\generalCompressibility}{\ensuremath{z}\xspace}
\newcommand{\compressibilityFactor}[1]{\ensuremath{\generalCompressibility_{#1}}\xspace}

\newcommand{\gravitationalAcceleration}{\ensuremath{g}\xspace}

\newcommand{\pipeLength}[1]{\ensuremath{L_{#1}}\xspace}
\newcommand{\pipeArea}[1]{\ensuremath{A_{#1}}\xspace}
\newcommand{\pipeDiameter}[1]{\ensuremath{D_{#1}}\xspace}
\newcommand{\pipeSlope}[1]{\ensuremath{s_{#1}}\xspace}
\newcommand{\generalFriction}{\ensuremath{\lambda}\xspace}
\newcommand{\frictionFactor}[1]{\ensuremath{\generalFriction_{#1}}\xspace}

\newcommand{\generator}{\ensuremath{\mathbb{G}_{\theta_{1}}}\xspace}
\newcommand{\discriminator}{\ensuremath{\mathbb{D}_{\theta_{2}}}\xspace}
\newcommand{\neuralnetwork}{\ensuremath{\mathbb{N}_{\{\theta_{1},\theta_{2}\}}}\xspace}

\newcommand{\mip}{\ensuremath{\mathbb{P}_{\pi}}\xspace}
\newcommand{\optmip}{\ensuremath{\mathbb{P}_{\pi}^{z_{1}^{*}}}\xspace}
\newcommand{\fixedmip}{\ensuremath{\mathbb{P}_{\pi}^{z_{1}}}\xspace}
\newcommand{\fixedapproxmip}{\ensuremath{\mathbb{P}_{\pi}^{\hat{z_{1}}}}\xspace}

\newcommand{\mipfeaturinggenerator}{\ensuremath{\mathbb{P}_{\pi}^{\generator(\pi)}}\xspace}

\newcommand{\LHS}{\ensuremath{A_{\pi}}\xspace}
\newcommand{\RHS}{\ensuremath{b_{\pi}}\xspace}
\newcommand{\binvars}{\ensuremath{z_{1}}\xspace}
\newcommand{\predbinvars}{\ensuremath{\hat{z_{1}}}\xspace}
\newcommand{\optbinvars}{\ensuremath{z_{1}^{*}}\xspace}
\newcommand{\instance}{\ensuremath{\pi}\xspace}
\newcommand{\instances}{\ensuremath{\Pi}\xspace}

\newcommand{\objmip}{\ensuremath{f(\mip)}\xspace}
\newcommand{\objfixedmip}{\ensuremath{f(\fixedmip)}\xspace}
\newcommand{\objfixedapproxmip}{\ensuremath{f(\fixedapproxmip)}\xspace}

\newcommand{\predobjfixedmip}{\ensuremath{\hat{f}(\fixedmip)}\xspace}
\newcommand{\predobjfixedapproxmip}{\ensuremath{\hat{f}(\fixedapproxmip)}\xspace}

\newcommand{\predobjmip}{\ensuremath{\hat{f}(\optmip)}\xspace}


%% file: sections/Abstract.tex
A decision support system relies on frequent re-solving of similar problem instances. While the general structure remains the same in corresponding applications, the input parameters are updated on a regular basis. We propose a generative neural network design for learning integer decision variables of mixed-integer linear programming (MILP) formulations of these problems. We utilise a deep neural network discriminator and a MILP solver as our oracle to train our generative neural network. In this article, we present the results of our design applied to the transient gas optimisation problem. With the trained network we produce a feasible solution in 2.5s, use it as a warm-start solution, and thereby decrease global optimal solution solve time by 60.5\%.

%% file: sections/Introduction.tex

Mixed-Integer Linear Programming (MILP) is concerned with the modelling and solving of problems from discrete optimisation. These problems can represent real-world scenarios, where discrete decisions can be appropriately captured and modelled by the integer variables. In real-world scenarios a MILP model is rarely solved only once. More frequently, the same model is used with varying data to describe different instances of the same problem which are solved on a regular basis. This holds true in particular for decision support systems, which can utilise MILP to provide real-time optimal decisions on a continual basis, see \cite{belien2009decision} and \cite{ruiz2004decision} for examples in nurse scheduling and vehicle routing. The MILPs that these decision support systems solve have identical structure due to both their underlying application and cyclical nature, and thus often have similar optimal solutions. Our aim is to exploit this repetitive structure, and create generative neural networks that generate binary decision encodings for subsets of important variables. These encodings can then be used in a primal heuristic by solving the induced sub-problem following variable fixations. Additionally, the then result of the primal heuristic can be used in a warm-start context to help improve solver performance in a globally optimal context. We demonstrate the performance of our neural network (NN) design on the transient gas optimisation problem \cite{rios2015optimization}, specifically on real-world instances embedded in day-ahead decision support systems. 

The design of our framework is inspired by the recent development of Generative Adversarial Networks (GANs) \cite{goodfellow2016nips}. Our design consists of two NNs, a Generator and a Discriminator. The Generator is responsible for generating the binary decision values, while the Discriminator is tasked with predicting the optimal objective function value of the MILP induced by fixing these binary variables to their generated values. 

Our NN design and its application to transient gas-network MILP formulations is an attempt to integrate Machine Learning (ML) into the MILP solving process. This integration has recently received an increased focus \cite{tang2019reinforcement, bertsimas2019online, gasse2019exact}, which has been encouraged by the success of ML integration into other facets of combinatorial optimisation, see \cite{bengio2018machine} for a thorough overview. Our contribution to this intersection of two fields is as follows: We introduce a new generative NN design for learning integer variables of parametric MILPs, which interacts with the MILP directly during training. We also apply our design to a much more difficult and convoluted problem than traditionally seen in similar papers, namely the transient gas transportation problem. This paper is to the best our knowledge the first successful implementation of ML applied to discrete control in gas-transport.

%% file: sections/RelatedWork.tex
As mentioned in the introduction, the intersection of MILP and ML is currently an area of active and growing research. For a thorough overview of Deep Learning (DL), the relevant subset of ML used throughout this article, we refer readers to \cite{goodfellow2016deep}, and for MILP to \cite{achterberg2007constraint}. We will highlight previous research from this intersection that we believe is either tangential, or may have shared applications to that presented in this paper. Additionally, we will briefly detail the state-of-the-art in transient gas transport, and highlight why our design is of practical importance. It should be noted as-well, that there are recent research activities aiming at the reverse direction, with MILP applied to ML instead of the orientation we consider, see \cite{wong2017provable} for an interesting example. 

Firstly, we summarise applications of ML to adjacent areas of the MILP solving process. \cite{gasse2019exact} creates a method for encoding MILP structure in a bipartite graph representing variable-constraint relationships. This structure is the input to a Graph Convolutional Neural Network (GCNN), which imitates strong branching decisions. The strength of their results stem from intelligent network design and the generalisation of their GCNN to problems of a larger size, albeit with some generalisation loss. \cite{zarpellon2020parameterizing} take a different approach, and use a NN design that incorporates the branch-and-bound tree state directly. In doing so, they show that information contained in the global branch-and-bound tree state is an important factor in variable selection. Furthermore, they are one of the few publications to present techniques on heterogeneous instances. \cite{etheve2020reinforcement} show a successful implementation of reinforcement learning for variable selection. \cite{tang2019reinforcement} show preliminary results of how reinforcement learning can be used in cutting-plane selection. By restricting themselves exclusively to Gomory cuts, they are able to produce an agent capable of selecting better cuts than default solver settings for specific classes of problems. 

There exists a continuous trade-off between model exactness and complexity in the field of transient gas optimisation, and as such, there is no standard model for transient gas transportation problems. \cite{moritz2007mixed} presents a piece-wise linear MILP approach to the transient gas transportation problem, \cite{burlacu2019maximizing} a non-linear approach with a novel discretisation scheme, and \cite{hennings2020controlling} and \cite{hoppmann2019optimal} a linearised approach. For the purpose of our experiments, we use the model of \cite{hennings2020controlling}, which uses linearised equations and focuses on active element heavy subnetworks. The current research of ML in gas transport is still preliminary. \cite{pourfard2019design} use a dual NN design to perform online calculations of a compressors operating point to avoid re-solving the underlying model. The approach constraints itself to continuous variables and experimental results are presented for a gunbarrel type network. \cite{mohamadibaghmolaei2014assessing} present a NN combined with a genetic algorithm for learning the relationship between compressor speeds and the fuel consumption rate in the absence of complete data. 
More often ML has been used in fields closely related to gas transport, as in \cite{hanachi2018performance}, with ML used to track the degradation of compressor performance, and in \cite{PetkovicChenGamrathetal.2019} to forecast demand values at the boundaries of the network. For a more complete overview of the transient gas literature, we refer readers to \cite{rios2015optimization}.

Our Discriminator design, which predicts the optimal objective value of an induced sub-MILP, can be considered similar to \cite{baltean2019scoring} in what it predicts and similar to \cite{ferber2019mipaal} in how it works. In the first paper \cite{baltean2019scoring}, a neural network is used to predict the associated objective value improvements on cuts. This is a smaller scope than our prediction, but is still heavily concerned with the MILP formulation. In the second paper \cite{ferber2019mipaal}, a technique is developed that performs backward passes directly through a MILP. It does this by solving MILPs exclusively with cutting planes, and then receiving gradient information from the KKT conditions of the final linear program. This application of a neural network, which produces input to the MILP, is very similar to our design. The differences arise in that we rely on a NN Discriminator to appropriately distribute the loss instead of solving a MILP directly, and that we generate variable values instead of parameter values with our Generator.

While our discriminator design is heavily inspired from GANs \cite{goodfellow2016nips}, it is also similar to actor-critic algorithms, see \cite{pfau2016connecting}. These algorithms have shown success for variable generation in MILP, and are notably different in that they sample from a generated distribution for down-stream decisions instead of always taking the decision with highest probability. Recently, \cite{chen2020actor} generated a series of coordinates for a set of UAVs using an actor-critic based algorithm, where these coordinates were continuous variables in a MINLP formulation. The independence of separable sub-problems and the easily realisable value function within their formulation resulted in a natural Markov Decision Process interpretation. For a better comparison on the similarities between actor-critic algorithms and GANs, we refer readers to \cite{pfau2016connecting}.

Finally, we summarise existing research that also deals with the generation of decision variable values for MIPs. \cite{bertsimas2018voice, bertsimas2019online} attempt to learn optimal solutions of parametric MILPs and MIQPs, which involves both outputting all integer decision variable values and the active set of constraints. They mainly use Optimal Classification Trees in \cite{bertsimas2018voice} and NNs in \cite{bertsimas2019online}. Their aim is tailored towards smaller problems classes, where speed is an absolute priority and parameter value changes are limited. \cite{masti2019learning} learn binary warm start decisions for MIQPs. They use NNs with a loss function that combines binary cross entropy and a penalty for infeasibility. Their goal of a primal heuristic is similar to ours, and while their design is much simpler, it has been shown to work effectively on very small problems.  Our improvement over this design is our non-reliance on labelled optimal solutions which are needed for binary cross entropy. \cite{ding2019optimal} present a GCNN design which is an extension of \cite{gasse2019exact}, and use it to generate binary decision variable values. Their contributions are a tripartite graph encoding of MILP instances, and the inclusion of their aggregated generated values as branching decisions in the branch-and-bound tree, both in an exact approach and in an approximate approach with local branching \cite{fischetti2003local}. Very recently, \cite{nair2020solving} combined the branching approach of \cite{gasse2019exact} with a novel neural diving approach, in which integer variable values are generated. They use a GCNN both for generating branching decisions and integer variables values. Different to our generator-discriminator based approach, they generate values directly from a learned distribution, which is based on an energy function that incorporates resulting objective values.

%% file: sections/NeuralNetworkDesign.tex
We begin by formally defining both a MILP and a NN. Our definition of a MILP is an extension of more traditional formulations, see \cite{achterberg2007constraint}, but still encapsulates general instances.
\begin{definition}\label{def:mip}
Let $\pi \in \mathbb{R}^p$ be a vector of problem defining parameters. We call the following a MILP parameterised by $\pi$.
\begin{align}
    \begin{split}
    \mip := \quad \text{min} & \quad  c_{1}^\mathsf{T}x_{1} + c_{2}^\mathsf{T}x_{2} + c_{3}^\mathsf{T}z_{1} + c_{4}^\mathsf{T}z_{2} \\
    \text{s.t} & \quad \LHS \begin{bmatrix} x_{1} \\ x_{2} \\ z_{1} \\ z_{2} \end{bmatrix} \leq \RHS \\
    & \quad c_{k} \in \mathbb{R}^{n_{k}}, k \in \{1,2,3,4\}, A_{\pi} \in \mathbb{R}^{m \times n}, b_{\pi} \in \mathbb{R}^{m} \\
    & \quad x_{1} \in \mathbb{R}^{n_{1}}, x_{2} \in \mathbb{R}^{n_{2}}, z_{1} \in \mathbb{Z}^{n_{3}}, z_{2} \in \mathbb{Z}^{n_{4}} \label{eq:MIP}
    \end{split}
\end{align}
Furthermore let $\Sigma \subset \mathbb{R}^p$ be a set of valid problem defining parameters. We then call $\{\mip \vert \pi \in \Sigma\}$ a problem class for $\Sigma$. 
\end{definition}
Note that the explicit parameter space $\Sigma$ is usually unknown, but we assume in the following to have access to a random variable $\Pi$ that samples from $\Sigma$. In addition, note that $c, n_{1}, n_{2}, n_{2}$, and $n_{4}$ are not parameterised by $\pi$, and as such the objective function and variable dimensions do not change between scenarios.

\begin{definition}\label{def:nn}
A $k$ layer NN $N_{\theta}$ is given by the following:
\begin{align}
\begin{split}
    N_{\theta} &: \mathbb{R}^{|a_{1}|} \xrightarrow{} \mathbb{R}^{|a_{k+1}|} \\
    h_{i} &: \mathbb{R}^{|a_{i}|} \xrightarrow{} \mathbb{R}^{|a_{i}|}, \quad \forall i \in \{2,...,k+1\} \\
    a_{i+1} &= h_{i+1}(W_{i}a_{i} + b_{i}),  \quad \forall i \in \{1,...,k\} 
\end{split}
\end{align}
Here $\theta$ fully describes all weights $(W)$ and biases $(b)$ of the network. $h_{i}$'s are called activation functions and are non-linear element-wise functions.
\end{definition}

\begin{figure}
    \centering
    \includegraphics[scale=0.5]{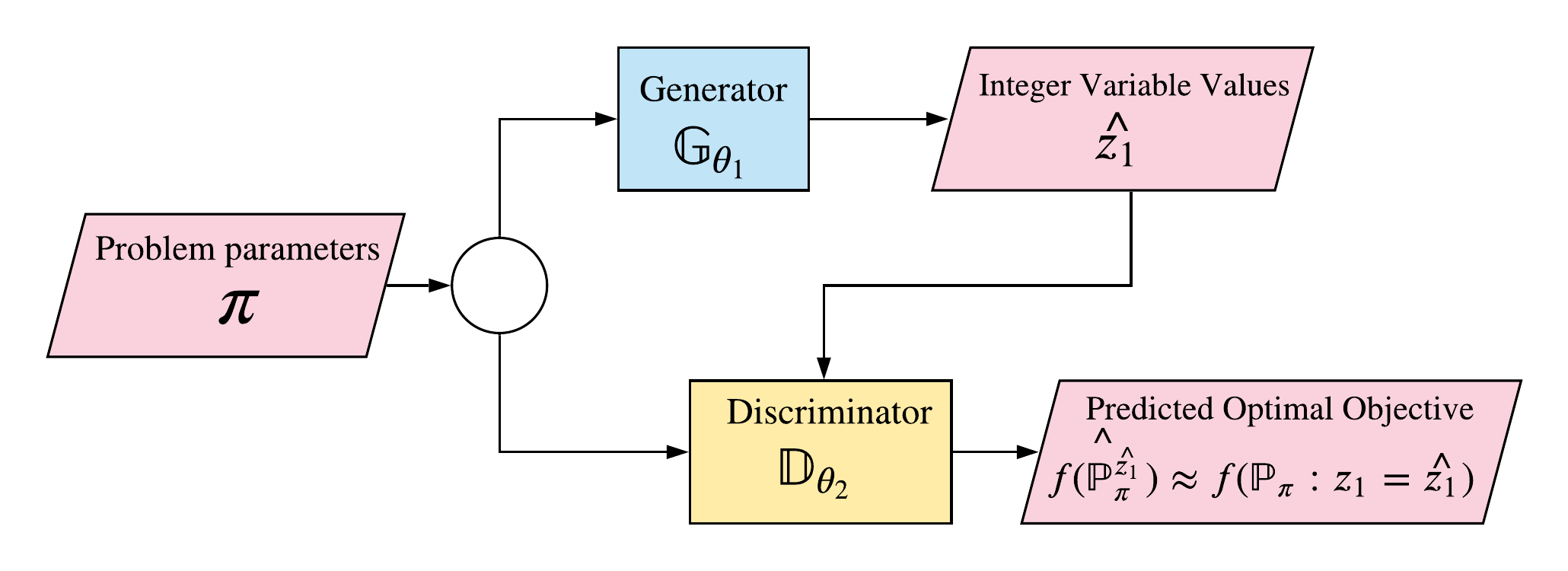}
    \caption{The general design of \neuralnetwork}
    \label{fig:forward_pass}
\end{figure}

An outline of our framework is depicted in Figure \ref{fig:forward_pass}. The Generator \generator is a NN that takes as input $\pi$. \generator outputs values for the variables \binvars, which we denote by \predbinvars. These variable values \predbinvars alongside $\pi$ are then input into another NN, namely the Discriminator \discriminator. \discriminator finally outputs a prediction of the optimal objective function value of \mip with values of \binvars fixed to \predbinvars, namely \predobjfixedapproxmip. More formally this is:

\begin{definition}
The generator \generator and discriminator \discriminator are both NNs defined by the following:
\begin{align}
    \begin{split}
        \generator &: \mathbb{R}^{p} \xrightarrow{} \mathbb{Z}^{n_{3}} \\
        \discriminator &: \mathbb{R}^{p} \times \mathbb{Z}^{n_{3}} \xrightarrow{} \mathbb{R}
    \end{split}
\end{align}
Furthermore, a forward pass of both \generator and \discriminator is defined as follows:
\begin{align}
    \hat{z_{1}} &= \generator(\pi) \label{eq:generator_defn}\\
    \predobjfixedapproxmip &= \discriminator(\hat{z_{1}},\pi) \label{eq:discriminator_defn} 
\end{align}
The hat notation is used to denote quantities that were approximated by a NN, and \objmip refers to the optimal objective function value of \mip. We use superscript notation to create the following instances:
\begin{align}
    \fixedapproxmip = \mip \quad \text{s.t} \quad z_{1} = \hat{z_{1}}
\end{align}
Note that the values of $\hat{z_{1}}$ must be appropriately rounded when explicitly solving \fixedapproxmip s.t they are feasible w.r.t. their integer constraints. As such, it is a slight abuse notation to claim that $\generator(\pi)$ lies in $\mathbb{Z}^{n_{3}}$
\end{definition}

The goal of this framework is to produce good initial solution values for \binvars, which lead to an induced sub-MILP, \fixedmip, whose optimal solution is a good feasible solution to the original problem.
Further, the idea is to use this feasible solution as a first incumbent for warm-starting \mip. To ensure feasibility for all choices of $\binvars$, we divide the continuous variables into two sets, $x_{1}$ and $x_{2}$, as seen in Definition \ref{def:mip}. The variables $x_{2}$ are potential slack variables to ensure that all generated decisions result in feasible $\fixedapproxmip$ instances. Penalising these slacks in the objective then feeds in naturally to our design, where $\generator$ aims to minimise the induced optimal objectives. For the purpose of our application it should be noted that $z_{1}$ and $z_{2}$ are binary variables instead of integer. Next we describe the design of \generator and \discriminator. 

\subsection{Generator and Discriminator Design}
\begin{figure}[h]
    \centering
    \includegraphics[width=0.95\linewidth]{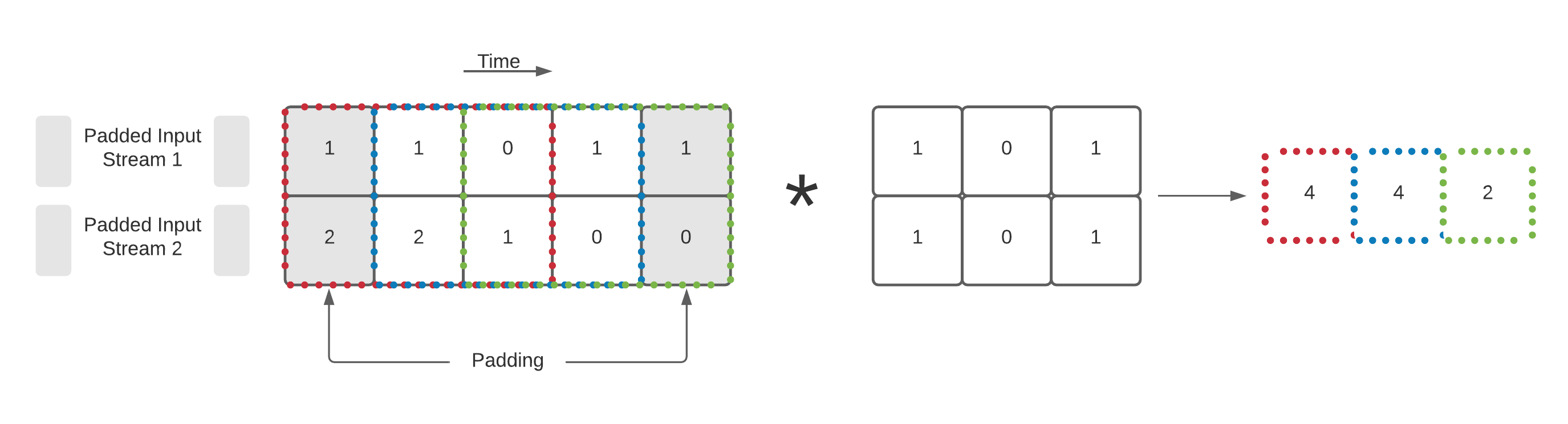}
    \caption{Method of merging two 1-D input streams}
    \label{fig:2d1d_convolution}
\end{figure}

\begin{figure}
    \centering
    \includegraphics[width=0.45\textwidth, height=0.5\textwidth]{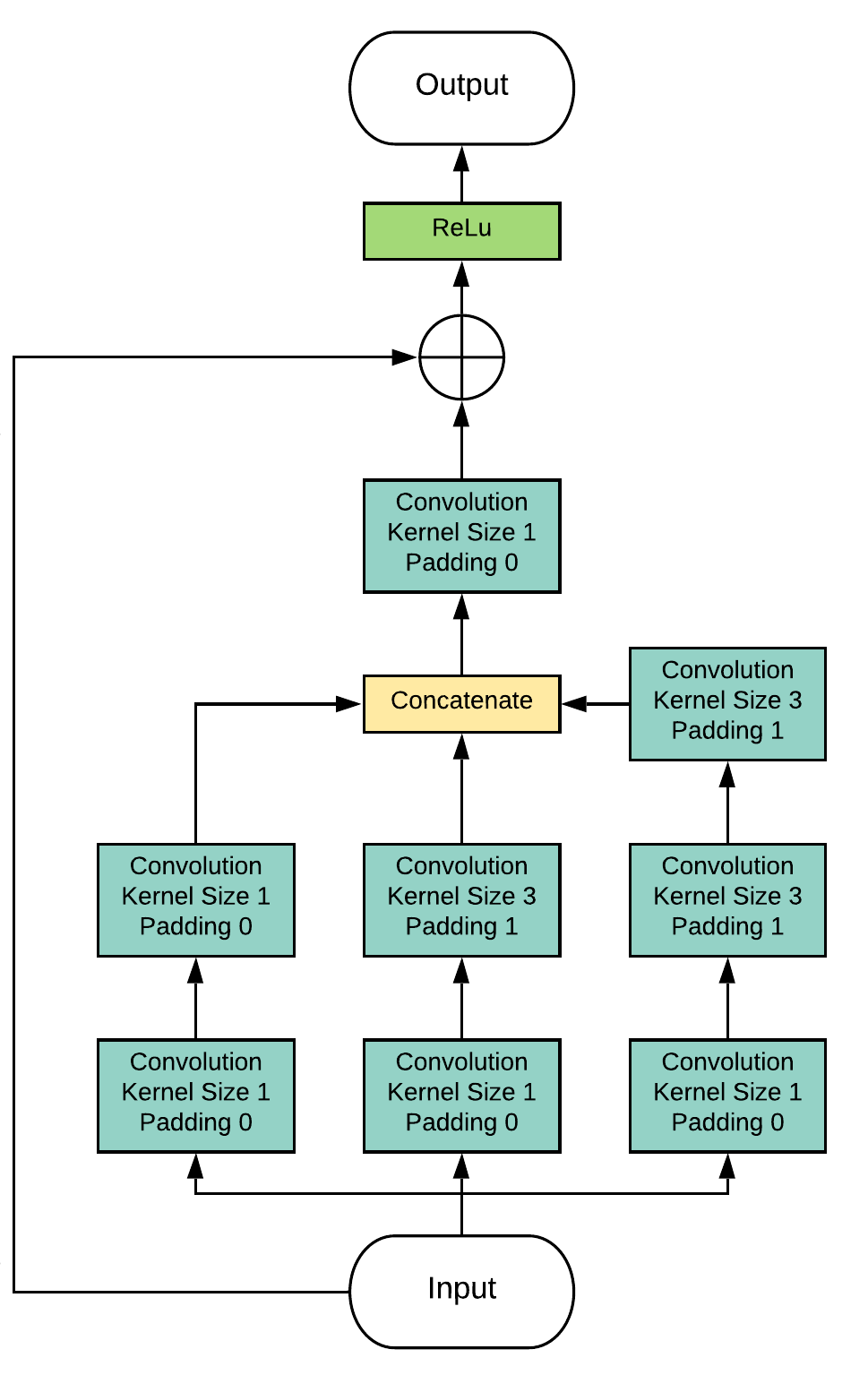}
    \caption{1-D Resnet-v2 Block Design}
    \label{fig:block_design}
\end{figure}

\generator and \discriminator are NNs whose structure is inspired by 
\cite{goodfellow2016nips}, as well as both inception blocks and residual NNs, which have greatly increased large scale model performance \cite{szegedy2017inception}.
We use the block design from Resnet-v2 \cite{szegedy2017inception}, see Figure \ref{fig:block_design}, albeit with slight modifications for the case of transient gas-network optimisation. Namely, we primarily use 1-D convolutions with that dimension being time. Additionally, we separate initial input streams by their characteristics, and when joining two streams, use 2-D convolutions, where the second dimension is of size 2 and quickly becomes one dimensional again. See Figure \ref{fig:2d1d_convolution} for an example of this process. 
The final layer of \generator contains a softmax activation function with temperature. As the softmax temperature increases, this activation function's output approaches a one-hot vector encoding. The final layer of \discriminator contains a softplus activation function. All other intermediate layers of \neuralnetwork use the ReLU activation function.
We refer readers to \cite{goodfellow2016deep} for a thorough overview of deep learning, and to Figure \ref{fig:graph_neural_network} in Appendix \ref{sec:appendix} for our complete design. 

For a vector $x=(x_{1}, \cdots, x_{n})$, the Softmax function with temperature $T \in \mathbb{R}$ \eqref{eq:softmax}, ReLu function \eqref{eq:relu}, and Softplus function with parameter $\beta \in \mathbb{R}$ \eqref{eq:softplus} are:
\begin{align}
    \sigma_{1}(x,T) &:= \frac{\exp(Tx_{i})}{\sum_{j=1}^{n}\exp(Tx_{j})} \label{eq:softmax} \\
    \sigma_{2}(x_{i}) &:= \max(0,x_{i}) \label{eq:relu} \\
    \sigma_{3}(x_{i},\beta) &:= \frac{1}{\beta}\log(1 + \exp(\beta x_{i})) \label{eq:softplus}
\end{align}

We can compose \generator with \discriminator, as in Figure
\ref{fig:forward_pass}, so that the combined resulting NN is defined as:
\begin{align}
    \neuralnetwork(\pi) := \discriminator(\generator(\pi),\pi)
\end{align} 

\subsection{Interpretations}

In a similar manner to GANs and actor-critic algorithms, see \cite{pfau2016connecting}, the design of \neuralnetwork has a bi-level optimisation interpretation, see \cite{dempe2002foundations} for an overview of bi-level optimisation. Here we list the explicit objectives of both \generator and \discriminator, and how their loss functions represent these objectives. 

The objective of \discriminator is to predict \objfixedapproxmip, the optimal induced objective values of \fixedapproxmip. Its loss function is thus:
\begin{align}
    L(\theta_{2}, \pi) :=  \big| \discriminator(\generator(\pi), \pi) - f(\mipfeaturinggenerator) \big| \label{eq:discriminator_approx_loss}
\end{align}

The objective of \generator is to minimise the induced prediction of \discriminator. Its loss function is thus:
\begin{align}
    L'(\theta_{1}, \pi) :=  \discriminator(\generator(\pi), \pi) \label{eq:generator_loss}
\end{align}

The corresponding bi-level optimisation problem can then be viewed as:

\begin{align}
    \begin{split}
    &\min\limits_{\theta_{1}} \quad \mathbb{E}_{\pi \sim \Pi} [ \discriminator(\generator(\pi), \pi) ] \\
    & \text{s.t} \quad \min\limits_{\theta_{2}} \quad \mathbb{E}_{\pi \sim \Pi} [ \discriminator(\generator(\pi), \pi) - f(\mipfeaturinggenerator) ]
    \end{split}
    \label{eq:bi-level}
\end{align}

\subsection{Training Method}
For effective training of \generator, a capable \discriminator is needed. We therefore pre-train \discriminator. The following loss function, which replaces $\generator(\instance)$ with prior generated \binvars values in \eqref{eq:discriminator_approx_loss}, is used for this pre-training:

\begin{align}
    L''(\theta_{2}, \pi) := \big| \discriminator(z_{1}, \pi) - \objfixedmip \big|
    \label{eq:discriminator_loss}
\end{align}

However, performing this initial training requires generating instances of \fixedmip. Here we do supervised training in an offline manner on prior generated data. 

After the initial training of \discriminator, we train \generator as a part of \neuralnetwork, using samples $\pi \in \Pi$, the loss function \eqref{eq:generator_loss}, and fixed $\theta_{2}$. The issue of \generator outputting continuous values for $\hat{z_{1}}$ is overcome by the final layer's activation function of \generator. The softmax with temperature \eqref{eq:softmax} ensures that adequate gradient information still exists to update $\theta_{1}$, and that the results are near binary. When using these results to explicitly solve \fixedapproxmip, we round our result to a one-hot vector encoding along the appropriate dimension. 

After the completion of both initial training, we alternately train both NN's using updated loss functions in the following way:

\begin{itemize}
    \item \discriminator training:
    \begin{itemize}
        \item As in the initial training, using loss function \eqref{eq:discriminator_loss}.
        \item In an online fashion, using predictions from \generator and loss function \eqref{eq:discriminator_approx_loss}.
    \end{itemize}
    \item \generator training:
    \begin{itemize}
        \item As explained above with loss function \eqref{eq:generator_loss}.
    \end{itemize}
\end{itemize}

Our design allows the loss to be back-propagated through \discriminator and distributed to the individual nodes of the final layer of \generator, i.e., that representing \binvars. This is largely different to other methods, many of which rely on using binary cross entropy loss against optimal solutions of $\mip$. Our advantage over these is that the contribution to the objective function we are trying to minimise of each variable decision in $z_{1}$ can be calculated. This has an added benefit of generated suboptimal solutions being much more likely to be near-optimal, as they are trained in a manner to minimise the objective rather than copy previously observed optimal solutions. 

For our application, transient gas network optimisation, methods for sampling instances currently do not exist. In fact, even gathering data is notoriously difficult, see  \cite{kunz2017electricity} and \cite{yueksel2020lessons}. For this reason, we introduce a new method for generating training data in section \ref{sec:Computational Experiments}.

%% file: sections/GasTransportModel.tex
To evaluate the performance of our approach, we test our framework on the transient gas optimisation problem, see \cite{rios2015optimization} for an overview of the problem and associated literature. This problem is difficult to solve as it combines a transient flow problem with complex combinatorics representing switching decisions. The natural modelling of transient gas networks as time-expanded networks lends itself well to our framework however, due to the static underlying network and repeated constraints at each time-step.

We use the description of transient gas networks by \cite{hennings2020controlling}. The advantages of this description for our framework is a natural separation of \binvars variables, which induce feasible \fixedmip for all choices due to the existence of slack variables in the description. These slack variables are then represented by $x_{2}$ in Definition \ref{def:mip}. The gas network is modelled as a directed graph \generalGraph where \setArcs is the set of arcs representing network elements, e.g. pipes, and the nodes \setVertices represent junctions between adjacent elements. Every arc $a \in \setArcs$ models a specific element with $\setArcsAllStationModel$, i.e., pipes, valves, resistors, regulators, and compressors. Additionally, the node set \setVertices contains multiple element types, with \setVerticesAll partitioned into boundary and inner nodes respectively. The boundary nodes represent the sources and sinks of the flow network. Thus, flow and pressure forecasts are given for each $v \in \setBoundaryNodes$.

It should be noted that this description focuses on \textit{network stations}, the beating hearts of gas networks. Network stations are commonly located at the intersections of major pipelines and contain nearly all elements, which can be used to control the gas flow. Next, we briefly explain the most important constraints from the model of \cite{hennings2020controlling}, particularly those which we exploit with our approach. For a full definition of the MILP, please see \cite{hennings2020controlling}.

As we optimise a transient problem, we deal with a time horizon, namely $\setTimestepsAll$. We aim to calculate a network state for each $t \in \setTimestepsNoZeroAll$, i.e. control decisions for all future time steps. As such, the initial gas network state at time 0 contains a complete description of that time step and is immutable. On the other hand all future time steps contain, before optimising, only forecasted pressure and flow values at $\setBoundaryNodes$. We denote \granularity{t} as the time difference in seconds from time step 0.

\subsection{Pipe Equations}
Pipes constitute the majority of elements in any gas transmission network. The dynamics of flow through pipes are governed by the Euler Equations, a set of nonlinear hyperbolic partial differential equations, see \cite{Osi1996}. We consider the isothermal case and discretise as in \cite{hennings2018benefits}. Consider the pipe $a=(u,v)$, $a \in \setPipes$, where $u, v \in \setVertices$ are the two incident nodes. We attach a flow-in $\varPipeFlow{u}{a}{t}$ and flow-out $\varPipeFlow{v}{a}{t}$ variable to each pipe. Additionally, each incident node has an attached pressure variable, namely  $(\varPressure{u}{t})$ and  $(\varPressure{v}{t})$. Moreover, these flow-in, flow-out, and pressure values also appear for each time step.
$R_{s}$, $z_{a}$, and \generalGasTemperature are assumed to be constant, and \pipeDiameter{a}, \pipeLength{a}, \pipeSlope{a}, \pipeArea{a}, \gravitationalAcceleration, and \frictionFactor{a} are themselves constant. The above constant assumptions are quite common in practice \cite{rios2015optimization}. 
It is only after setting the velocity of gas within each individual pipe, $\paramAbsoluteVelocity{w}{a}$ to be constant that all non-linearities are removed however. We do this via a method developed in \cite{hennings2018benefits} and seen in \cite{fang2017dynamic}. The resulting pipe equations are:

\begin{align}
    \varPressure{u}{t_2} + \varPressure{v}{t_2} - \varPressure{u}{t_1} - \varPressure{v}{t_1}
    + \frac{2\specificGasConstant\generalGasTemperature\compressibilityFactor{a}(\granularity{t_2} - \granularity{t_1})}{\pipeLength{a}\pipeArea{a}}
      \left(\varPipeFlow{v}{a}{t_2} - \varPipeFlow{u}{a}{t_2}\right) &=0 \label{eq:pipes_constVelo_continuity}\\
    \varPressure{v}{t_2} - \varPressure{u}{t_2}
  + \frac{\frictionFactor{a}\pipeLength{a}}{4\pipeDiameter{a}\pipeArea{a}}
    \left(\paramAbsoluteVelocity{u}{a}\varPipeFlow{u}{a}{t_2} + \paramAbsoluteVelocity{v}{a}\varPipeFlow{v}{a}{t_2}\right) & \notag\\
  + \frac{\gravitationalAcceleration\pipeSlope{a}\pipeLength{a}}{2\specificGasConstant\generalGasTemperature\compressibilityFactor{a}}
    \left(\varPressure{u}{t_2} + \varPressure{v}{t_2}\right) &=0 \label{eq:pipes_constVelo_momentum}
\end{align}

As nodes represent junctions between network elements and thus have no volume in which to store any gas, the flow conservation constraints \eqref{eq:flowConservation_boundaryNodes} \eqref{eq:flowConservation_innerNodes} are required. In the below equations, $\varInflowValue{v}{t}$ represents the inflow resp. outflow of entry and exit nodes in the network at time $t \in \setTimesteps$. Note that network elements that aren't pipes have only one associated flow variable, instead of the in-out flow exhibited by pipes. This is due to them having no volume, and as such no ability to store gas over time, i.e. line-pack. 

\begin{align}
     & \sum_{(u,w)=a\in\setPipes} \varPipeFlow{w}{a}{t} - \sum_{(w,v)=a\in\setPipes} \varPipeFlow{w}{a}{t} \notag\\
    +& \sum_{(u,w)=a\in\setArcs\setminus\setPipes} \varArcFlow{a}{t} - \sum_{(w,v)=a\in\setArcs\setminus\setPipes} \varArcFlow{a}{t} + \varInflowValue{w}{t} = 0 \qquad \forall w\in\setBoundaryNodes \label{eq:flowConservation_boundaryNodes}\\
     & \sum_{(u,w)=a\in\setPipes} \varPipeFlow{w}{a}{t} - \sum_{(w,v)=a\in\setPipes} \varPipeFlow{w}{a}{t} \notag\\
    +& \sum_{(u,w)=a\in\setArcs\setminus\setPipes} \varArcFlow{a}{t} - \sum_{(w,v)=a\in\setArcs\setminus\setPipes} \varArcFlow{a}{t} = 0 \qquad \forall w\in\setInnerNodes \label{eq:flowConservation_innerNodes}
\end{align}

\subsection{Operation Modes}
\textit{Operation modes} represent binary decisions in our gas network. We identify the corresponding binary variables with the \binvars variables from our MILP formulation \eqref{eq:MIP}.
Let $\setNSModes$ represent the set of operation modes, and $\varNSmode{o}{t}$ the associated variables. Operation Modes are very important in our modelling context as they describe every allowable combination of discrete decisions associated with \textit{valves} and \textit{compressors}. 

\subsubsection{Compressors}
Compressors are typically set up as a compressor station consisting of multiple compressor units, which represent the union of one single compressor machine and its associated drive. These compressor units are dynamically switched on or off and used in different sequences to meet the current needs in terms of compression ratios and flow rates. Out of the theoretically possible arrangements of compressor units, the set of technically feasible arrangements are known as the \textit{configurations} of a compressor station.

Selecting an operation mode results in fixed configurations for all compressor stations. The binary variables associated with a compressor station $a = (u,v) \in \setCompressorStations$ at time $t \in \setTimesteps$ are $\varModeBypass{a}{t}$ (bypass), $\varModeClosed{a}{t}$ (closed), and $\varModeConfiguration{c}{a}{t}$ $\forall c \in \setCompressorConfigurations{a}$ (active). $\setCompressorConfigurations{a}$ denotes the set of configurations associated to compressor station $a$ available in active mode, where the configuration's operating range is a polytope in space $(\varPressure{u}{t}, \varPressure{v}{t}, \varPipeFlow{u}{a}{t})$. The polytope of configuration $c$ is represented by the intersection of half-spaces, $\setConfigurationFacets{c} = \{(\alpha_0,\alpha_1,\alpha_2,\alpha_3) \in \mathbb{R}^4\}$.

\begin{align}
    1 &= \sum_{c\in\setCompressorConfigurations{a}} \varModeConfiguration{c}{a}{t} + \varModeBypass{a}{t} + \varModeClosed{a}{t} \label{eq:compressorStation_OneModeOrConfig} \\
    \begin{split}
    \alpha_0\varConfigPressureL{c}{a}{t} + \alpha_1\varConfigPressureR{c}{a}{t} + \alpha_2\varConfigFlow{c}{a}{t} + \alpha_3\varModeConfiguration{c}{a}{t} &\leq 0 \\
    & \forall (\alpha_0,\alpha_1,\alpha_2,\alpha_3) \in\setConfigurationFacets{c} \quad \forall c\in\setCompressorConfigurations{a} \label{eq:compressorStation_cfg_facets}
    \end{split}
\end{align}

Note that the variables in \eqref{eq:compressorStation_cfg_facets} have an extra subscript and superscript compared to those in \eqref{eq:pipes_constVelo_continuity} and \eqref{eq:pipes_constVelo_momentum}. This is due to our use of the convex-hull reformulation, see \cite{Bal2018}. The additional subscript refers to the configuration in question, and the superscript the mode, with the pressure variables having an additional node identifier. It should also be noted that the continuous variables attached to a compressor station are not fixed by a choice in operation mode or configuration, but rather the operation mode restricts the variables to some polytope. 

\subsubsection{Valves}
\textit{Valves} decide the allowable paths through a network, and can separate areas, decoupling their pressure levels. They are modelled as an arc $a=(u,v)$, whose discrete decisions can be decided by an operation mode choice. Valves have two modes, namely open and closed. When a valve is open, similar to a compressor station in bypass, flow is unrestricted and there exists no pressure difference between the valves start and endpoints. Alternatively in the closed mode, a valve allows no flow to pass, and decouples the pressure of the start- and endpoints of the arc. The variable $\varModeOpen{a}{t}$ represents a valve being open with value 1 and closed with value 0. The general notation \lb{x} and \ub{x} refer to lower and upper bounds of a variable $x$. The constraints describing valves are then as follows:

\begin{align}
    \varPressure{u}{t} - \varPressure{v}{t} &\leq ( 1 - \varModeOpen{a}{t} )(\paramPressureUB{u}{t}-\paramPressureLB{v}{t}) \label{eq:valves_first}\\
    \varPressure{u}{t} - \varPressure{v}{t} &\geq ( 1 - \varModeOpen{a}{t} )(\paramPressureLB{u}{t}-\paramPressureUB{v}{t}) \\
    \varArcFlow{a}{t} &\leq ( \varModeOpen{a}{t} )\paramArcFlowUB{a}{t} \\
    \varArcFlow{a}{t} &\geq ( \varModeOpen{a}{t} )\paramArcFlowLB{a}{t}. \label{eq:valves_last} 
\end{align}

\subsubsection{Valid Operation Modes}
As mentioned earlier, not all combinations of compressor station configurations and valve states are possible. We thus define a mapping $M(o,a)$ from operation mode $o \in \setNSModes$ to the discrete states of all $a \in \setValves \cup \setCompressorStations$

\begin{align*}
    M(o,a) := m &\text{ where }m\text{ is the mode or configuration of arc }a\\
    & \text{ in operation mode }o \quad \forall o\in\setNSModes\quad\forall a\in\setValves\cup\setCompressorStations\\
    \text{with}\qquad & m\in \{\text{op}, \text{cl}\} \text{ if }a\in\setValves \\
                      & m\in \{\text{by}, \text{cl}\} \cup \setCompressorConfigurations{a} \text{ if }a\in\setCompressorStations
\end{align*}

Using this mapping we can then define a set of constraints for all valid combinations of compressor station and valve discrete states for each $t \in \setTimestepsNoZero$. The variable \varNSmode{o}{t}, $o \in \setNSModes$ $t \in \setTimestepsNoZero$, is a binary variable, where the value 1 represents the selection of $o$ at time step $t$. 

\begin{align}
    \sum_{o \in \setNSModes} \varNSmode{o}{t} &= 1 \label{eq:opMode_choice} \\
    \varModeOpen{a}{t}   &= \sum_{o\in\setNSModes : M(o,a)=\text{op}} \varNSmode{o}{t} \quad\forall a\in\setValves \label{eq:valve_opMode_coupling}\\
    \varModeBypass{a}{t} &= \sum_{o\in\setNSModes : M(o,a)=\text{by}} \varNSmode{o}{t} \quad\forall a\in\setCompressorStations \label{eq:cs_bypass_opMode_coupling} \\
    \varModeClosed{a}{t} &= \sum_{o\in\setNSModes : M(o,a)=\text{cl}} \varNSmode{o}{t} \quad\forall a\in\setCompressorStations \label{eq:cs_closed_opMode_coupling} \\
    \varModeConfiguration{c}{a}{t} &= \sum_{o\in\setNSModes : M(o,a)=c}   \varNSmode{o}{t} \quad\forall c\in\setCompressorConfigurations{a}\quad\forall a\in\setCompressorStations \label{eq:cs_cfg_opMode_coupling} \\
    \varNSmode{o}{t} &\in \{0,1\} \quad \forall o\in\setNSModes. \notag
\end{align}

\subsection{Flow Directions}
\textit{Flow Directions} define the sign of flow values over the boundary nodes of a network station. With regards to our MILP they are a further set of decision variables. We avoid generating these decisions with our deep learning framework as not all combinations of operation modes and flow directions are feasible. These variables thus exist as integer variables in \fixedmip, namely as a subset of $z_{2}$, see \eqref{eq:MIP}. They are few in number however due to the limited combinations after the operation modes are fixed.

\subsection{Boundary Nodes and Slack} \label{sec:Boundary Nodes and Slaack}
Boundary nodes, unlike inner nodes, have a prescribed flow and pressure values for all future time steps. For each boundary node $v \in \setBoundaryNodes$ and $t \in \setTimestepsNoZero$, we have \varSlackPressurePos{v}{t} and \varSlackPressureNeg{v}{t}, which capture the positive and negative difference between the prescribed 
and realised pressure. In addition to these pressure slack variables, we have the inflow slack variables \varSlackFlowPos{v}{t} and \varSlackFlowNeg{v}{t} which act in a similar manner but for inflow. The relationships between the slack values, prescribed values, and realised values can be modelled for each $v \in \setBoundaryNodes$ and $t \in \setTimestepsNoZero$ as:
\begin{align}
  \demandPressure{v}{t} &= \varPressure{v}{t} - \varSlackPressurePos{v}{t} + \varSlackPressureNeg{v}{t} \quad \forall v\in\setBoundaryNodes \label{eq:slack_pressure}\\
  \demandInflow{v}{t} &= \varInflowValue{v}{t} - \varSlackFlowPos{v}{t} + \varSlackFlowNeg{v}{t} \quad \forall v \in \setBoundaryNodes \label{eq:slack_flow} 
\end{align}

Note that unlike the model from \cite{hennings2020controlling}, we do not allow the inflow over a set of boundary nodes to be freely distributed according to which group they belong to. This is an important distinction, as each single node has a complete forecast.

\subsection{Initial State}
In addition to the forecast mentioned in subsection \ref{sec:Boundary Nodes and Slaack}, we also start our optimisation problem with an initial state. This initial state contains complete information of all discrete states and continuous values for all network elements at $t=0$. 

\subsection{Objective function}
The objective of our formulation is to both minimise slack usage, and changes in network operation. Specifically, it is a weighted sum of changes in the active element modes, changes in the continuous active points of operation, and the deviations from given pressure and flow demands. For the exact objective function we refer readers to \cite{hennings2020controlling}. 

%% file: sections/ComputationalExperiments.tex
In this section we propose an experimental design to determine the effectiveness of our neural network design approach. We outline how we generate synthetic training data, and show the exact architecture and training method we use for our neural network. Our final test set consists of 15 weeks of real-world data provided by our project partner OGE.

\subsection{Data Generation}
As mentioned previously, acquiring gas network data is notoriously difficult \cite{yueksel2020lessons, kunz2017electricity}. Perhaps because of this difficulty, there exists no standard method for generating valid states for a fixed gas network. Below we outline our methods for generating synthetic transient gas instances for training purposes, i.e. generating $\instance \in \instances$ and artificial \binvars values. For our application of transient gas instances, \instance is a tuple of a boundary forecast and an initial state.

\subsubsection{Boundary Forecast Generation}
\label{subsection:Boundary Forecast Generation}
We consider network stations as our gas network topology. They contain all heavy machinery and at most only short segments of large scale transport pipelines. As such, our gas networks cannot be used to store large amounts of gas. We thus aim to generate balanced demand scenarios, with the requirement described as follows: 

\begin{align}
    \sum_{v \in \setBoundaryNodes} \demandInflow{v}{t} = 0 \quad \forall t \in \setTimestepsNoZero \label{eq:forecast_1}
\end{align}

The distribution of gas demand scenarios is not well known. Hence we naively assume a uniform distribution, and using the largest absolute flow value found over any node and time step in our real-world data, create an interval as follows:

\begin{align}
\begin{split}
    M_{\text{q}} &= \max_{v \in \setBoundaryNodes, t \in \setTimestepsNoZero} | \demandInflow{v}{t} | \\
    \demandInflow{v}{t} &\in \interval{-1.05M_{\text{q}}}{1.05M_{\text{q}}}
\end{split}\label{eq:forecast_2}
\end{align}

In addition to the above, we require three MILP formulation specific requirements. The first is that the absolute difference between the flow values of a node is not too large for any adjacent time steps. Secondly, the sign of the generated flow values must match the attribute of the boundary node, i.e., entry (+), exit (-). Thirdly, the flow values do not differ too largely between boundary nodes of the same \textit{fence group} within the same time step. A fence group is denoted by $g\in \setNSFenceGroups$, and enforces the sign of all nodes in the group to be identical. These constraints are described below:

\begin{align}
    \begin{split}
        | \demandInflow{v}{t} - \demandInflow{v}{t-1} | &\leq 200 \quad \forall t \in \setTimestepsNoZero, \quad v \in \setBoundaryNodes \\
        \text{sign}(\demandInflow{v}{t}) &= 
        \begin{cases}
            1 \quad \text{if} \quad v \in \setEntries \\
            -1 \quad \text{if} \quad v \in \setExits
        \end{cases} \forall t \in \setTimestepsNoZero, \quad v \in \setBoundaryNodes \\
        | \demandInflow{v_{1}}{t} - \demandInflow{v_{2}}{t} | &\leq 200 \quad \forall t \in \setTimestepsNoZero, \quad v_{1},v_{2} \in g, \hspace{0.5em} g\in \setNSFenceGroups, \hspace{0.5em} v_{1},v_{2} \in \setBoundaryNodes
    \end{split} \label{eq:forecast_3}
\end{align}

To generate demand scenarios that satisfy constraints \eqref{eq:forecast_1} and \eqref{eq:forecast_2}, we use the method proposed in \cite{rubin1981bayesian}. Its original purpose was to generate samples from the Dirichlet distribution, but it can be used for a special case of the Dirichlet distribution that is equivalent to a uniform distribution over a simplex in 3-dimensions. Such a simplex is exactly described by \eqref{eq:forecast_1} and \eqref{eq:forecast_2} for each time step. Hence we can apply it for all time-steps and reject all samples that do not satisfy constraints \eqref{eq:forecast_3}. Note that this method is insufficient for network stations with more than three boundary nodes.  

In addition to flow demands, we require a pressure forecast for all boundary nodes. Our only requirements here is that the pressures between adjacent time steps for a single node not fluctuate heavily and that the bounds are respected. We create a bound on the range of pressure values by finding maximum and minimum values over all nodes and time steps in our test set. We once again assume our samples to be uniformly distributed and sample appropriately over \eqref{eq:forecast_4} with rejection of samples that do not respect constraint \eqref{eq:forecast_5}. Note that many scenarios generated by this approach are unlikely to happen in practice, as the pressure and flow profiles may not match. 

\begin{align}
    \begin{split}
        M_{\text{p}}^{+} &= \max_{v \in \setBoundaryNodes, t \in \setTimestepsNoZero} \demandPressure{v}{t} \quad \quad
        M_{\text{p}}^{-} = \min_{v \in \setBoundaryNodes, t \in \setTimestepsNoZero} \demandPressure{v}{t} \\
        \demandPressure{v}{t} &\in \interval{M_{\text{p}}^{-} - 0.05(M_{\text{p}}^{+} - M_{\text{p}}^{-})}{M_{\text{p}}^{+} + 0.05(M_{\text{p}}^{+} - M_{\text{p}}^{-})}
    \end{split} \label{eq:forecast_4} 
\end{align}
\begin{align}
    | \demandPressure{v}{t} - \demandPressure{v}{t-1} | &\leq 5 \quad \forall t \in \setTimestepsNoZero, \quad v \in \setBoundaryNodes \label{eq:forecast_5}
\end{align}

Combining the two procedures from above yields the artificial forecast data generation method described in Algorithm \ref{alg:boundary_data_generator}. 

\begin{algorithm}[H]
\SetAlgoLined
\KwResult{A forecast of pressure and flow values over the time horizon}
flow\_forecast = Sample simplex \eqref{eq:forecast_1}\eqref{eq:forecast_2} uniformly, rejecting via \eqref{eq:forecast_3} \;
pressure\_forecast = Sample \eqref{eq:forecast_4} uniformly, rejecting via \eqref{eq:forecast_5} \;
\Return (flow\_forecast, pressure\_forecast)
\caption{Boundary Value Forecast Generator}
\label{alg:boundary_data_generator}
\end{algorithm}

\subsubsection{Operation Mode Sequence Generation}
\label{subsection:Operation Mode Sequence Generation}
During offline training, \discriminator requires optimal solutions for a fixed $z_{1}$. In Algorithm \ref{alg:op_mode_generator} we outline a naive yet effective approach of generating reasonable $z_{1}$ values, i.e., operation mode sequences: 

\begin{algorithm}
\SetAlgoLined
\KwResult{An Operation Mode per time step}
 operation\_modes = $\lbrack$ $\rbrack$ \; 
    \For{$t = 1;\ t < | \setTimestepsNoZero |;\ t = t + 1$}{
        \uIf{$t == 1$}{
            new\_operation\_mode = rand($\setNSModes$) \;
        }
        \ElseIf{rand(0,1) $\geq$ 0.9}{
            new\_operation\_mode = rand($\setNSModes \setminus$ new\_operation\_mode) \;
        }
    operation\_modes.append(new\_operation\_mode) \;
    }
    \Return operation\_modes
 \caption{Operation Mode Sequence Generator}
 \label{alg:op_mode_generator}
\end{algorithm}

\subsubsection{Initial State Generation}
\label{subsection:Initial State Generation}
Many coefficients of $A_{\pi}$ are invariant due to static network topology. Many others however are found by substituting multiple parameters into an equation describing gas properties. This information is contained in the initial state, and we generate them similar to boundary forecasts:
\begin{align}
\begin{split}
c_{\text{state}} &\in \text{$\{$Temperature, Inflow Norm Density, Molar Mass,} \\
\text{Pseudo } & \text{Critical Temperature, Pseudo Critical Pressure$\}$} 
\end{split} \\
\begin{split}
    M_{\text{c}}^{+} &= \max_{\text{state} \in \text{initial states}}  c_{\text{state}} \quad \quad M_{\text{c}}^{-} = \min_{\text{state} \in \text{initial states}} c_{\text{state}} \\
    c_{\text{state}} &\in \interval{M_{\text{c}}^{-} - 0.05(M_{\text{c}}^{+} - M_{\text{c}}^{-})}{M_{\text{c}}^{+} + 0.05(M_{\text{c}}^{+} - M_{\text{c}}^{-})}
\end{split}\label{eq:generated_constants}
\end{align}

We now have the tools to generate synthetic initial states, see Algorithm \ref{alg:initial_state_generator}.

Algorithm \ref{alg:initial_state_generator} is designed to output varied and valid initial states w.r.t our MILP formulation. However, it comes with some drawbacks. Firstly, the underlying distribution of demand scenarios for both flow and pressure are probably not uniform nor conditionally independent. Moreover, the sampling range we use is significantly larger than that of our test set as we take single maximum and minimum values over all nodes. Secondly, the choice of operation modes that occur in reality is also not uniform. In reality, some operation modes occur with a much greater frequency than others. Our data is thus more dynamic than reality, and likely to contain operation mode choices that do match the demand scenarios. Finally, we rely on a MILP solver to generate new initial states in our final step. Hence we cannot rule out the possibility of a slight bias. One example would be the case of a repeated scenario, which has multiple optimal solutions, but the MILP solver always returns an identical solution. 

\begin{algorithm}[H]
\SetAlgoLined
\SetKwInOut{Input}{Input}
\Input{Desired time-step distance $j \in [1, \cdots, k]$}
\KwResult{An initial state to the transient gas optimisation problem}
flow\_forecast, pressure\_forecast = Boundary Prognosis Generator() \footnote{See Algorithm \ref{alg:boundary_data_generator}} \;
gas\_constants = Sample
\eqref{eq:generated_constants} uniformly \;
initial\_state = Select random state from real-world data \;
\instance = (flow\_forecast, pressure\_forecast, gas\_constants, initial\_state) \footnote{Note that in general our \instance does not include gas\_constants. This is because the information is generally encoded in initial\_state. Our gas\_constants in this context are randomly generated however, and may not match the initial\_state. This does not affect solving as these values are simply taken as truths.} \;
\binvars =  Operation Mode Sequence Generator() \footnote{See Algorithm \ref{alg:op_mode_generator}} \;
\fixedmip = generate from \instance and \binvars \;
$($ state\_1, $\cdots$, state\_k $)$ = Optimal solution states from solving \fixedmip \;
\Return state\_j 
\caption{Initial State Generator}
\label{alg:initial_state_generator}
\end{algorithm}

\begin{algorithm}[H]
\SetAlgoLined
\SetKwInOut{Input}{Input}
\Input{num\_states, num\_scenarios, time\_step\_difference}
\KwResult{num\_scenarios many gas instances and their optimal solutions}
initial\_states = [] \;
\For{$i = 0;\ i < \text{num\_states};\ i = i + 1$}{
initial\_states.append(Initial State Generator(time\_step\_difference))\footnote{See Algorithm \ref{alg:initial_state_generator}}\;
}
forecasts = [] \;
\For{$i = 0;\ i < \text{num\_scenarios};\ i = i + 1$}{
flow\_forecast, pressure\_forecast = Boundary Prognosis Generator()\footnote{See Algorithm \ref{alg:boundary_data_generator}}\;
forecasts.append((flow\_forecast, pressure\_forecast)) \;
}
solve\_data = [] \;
\For{$i = 0;\ i < \text{num\_scenarios};\ i = i + 1$}{
\binvars = Operation Mode Sequence Generator() \footnote{See Algorithm \ref{alg:op_mode_generator}} \;
initial\_state = Uniformly select from initial\_states \;
\instance = (forecasts[i], initial\_state) \; 
\fixedmip = Create MILP from \instance and \binvars \;
solution = Solve \fixedmip \;
solve\_data.append((\binvars, \instance, solution)) \;
}
\Return solve\_data
\caption{Synthetic Gas Data Generator}
\label{alg:data_generator}
\end{algorithm}

In the case of initial state generation, we believe that further research needs to be performed. Our method is effective in the context of machine learning where we aim for a diverse set of data, but it is naive and incapable of ensuring that generated boundary scenarios are realistic. 

\subsubsection{Complete Transient Gas Instance Generation}
To train \discriminator and \generator, we need both the transient gas transportation scenario, and an optimal solution for it. Combining the generation methods for synthetic data in subsections \ref{subsection:Boundary Forecast Generation}, \ref{subsection:Operation Mode Sequence Generation},  \ref{subsection:Initial State Generation}, and the solving process of the created instances, we derive Algorithm \ref{alg:data_generator}. 

\subsection{Experimental Design}
We generated our initial training and validation sets offline. To do so we use Algorithm \ref{alg:data_generator} with inputs:
num\_states = $10^{4}$, num\_scenarios = $4\times 10^{6}$, and time\_step\_difference = 8.
This initial training data is exclusively used for training \discriminator, and is split into a training set of size $3.2\times 10^{6}$, a test set of $4\times 10^{5}$, and a validation set of $4\times 10^{5}$. 

The test set is checked against at every epoch, while the validation set is only referred to at the end of the initial training. Following this initial training, we begin to train \neuralnetwork as a whole, alternating between \generator and \discriminator. The exact algorithm is given in \ref{alg:train}, which references functions provided in Appendix \ref{sec:appendix}. For training, we used the Adam algorithm \cite{kingma2014adam} as our descent method. The associated parameters to this algorithm and a complete set of other training parameters are listed in Table \ref{tab:initial_discriminator}. In the case of a parameter being non-listed, the default value was used. The intention behind our training method is to ensure that \neuralnetwork receives no real-world data prior to its final evaluation. With this method we hope to show that synthetic data is sufficient for training purposes and that \neuralnetwork successfully generalises to additional data sets. However, we should note that Algorithm \ref{alg:initial_state_generator} does use real-world data as a starting point from which to generate artificial data.

\newpage

\input{sections/TrainingAlgo}

We consider the solution of \fixedapproxmip as a primal heuristic for the original problem \mip. Due to our usage of slack, i.e. the application of variables $x_{2}$, any valid solution for \fixedmip is a valid solution of \mip. We aim to incorporate \neuralnetwork in a global MIP context and do this by using a partial solution of \fixedapproxmip as a warm-start suggestion for \mip. The partial solution consists of \predbinvars, an additional set of binary variables called the flow directions, which are a subset of $z_{2}$ in \eqref{eq:MIP}, and \varPressure{v}{t} $ \forall v \in \setBoundaryNodes, t \in \setTimestepsNoZero$, which are a subset of $x_{1}$ in \eqref{eq:MIP}. Note that partial solutions are used as instances are numerically difficult. In doing so, we hope to generate valid solutions quickly, and speed up the global solution process. The primal heuristic and warm-start algorithm can be seen in Algorithms \ref{alg:primal_heuristic} and \ref{alg:warm_start} respectively. \\

\begin{algorithm}[H]
\SetAlgoLined
\SetKwInOut{Input}{Input}
\Input{\mip}
\predbinvars = \generator(\instance) \;
\fixedapproxmip = Create MILP from \instance and \predbinvars  \;
solution = Solve \fixedapproxmip \;
\Return solution \;
\KwResult{Optimal solution of \fixedapproxmip, primal solution of \mip.}
\caption{Primal Heuristic}
\label{alg:primal_heuristic}
\end{algorithm}

\begin{algorithm}[H]
\SetAlgoLined
\SetKwInOut{Input}{Input}
\Input{\mip}
primal\_solution = Primal Heuristic(\mip) \footnote{See Algorithm \ref{alg:primal_heuristic}} \;
optimum = Solve \mip with primal\_solution as a warm-start suggestion \;
\KwResult{Optimal solution of \mip}
\caption{Warm Start Algorithm}
\label{alg:warm_start}
\end{algorithm}

For our experiments we used PyTorch 1.4.0 \cite{paszke2019pytorch} as our ML modelling framework, Pyomo v5.5.1 \cite{hart2017pyomo, hart2011pyomo} as our MILP modelling framework, and Gurobi v9.02 \cite{gurobi} as our MILP solver. The MILP solver settings are available in Table \ref{tab:mip_params} in Appendix \ref{sec:appendix}. \neuralnetwork was trained on a machine running Ubuntu 18, with 384\,GB of RAM, composed of 2x \emph{Intel(R) Xeon(R) Gold 6132} running $@$ 2.60GHz, and 4x \emph{NVIDIA Tesla V100 GPU-NVTV100-16}. The final evaluation times were performed on a cluster using 4 cores and 16\,GB of RAM of a machine composed of 2x \emph{Intel Xeon CPU E5-2680} running $@$ 2.70\,GHz.

Our validation set for the final evaluation of \neuralnetwork consists of 15 weeks of live real-world data from our project partner OGE. Instances are on average 15 minutes apart for this period and total 9291. 

\begin{figure}[ht]
    \centering
    \includegraphics[scale=0.2, angle=0]{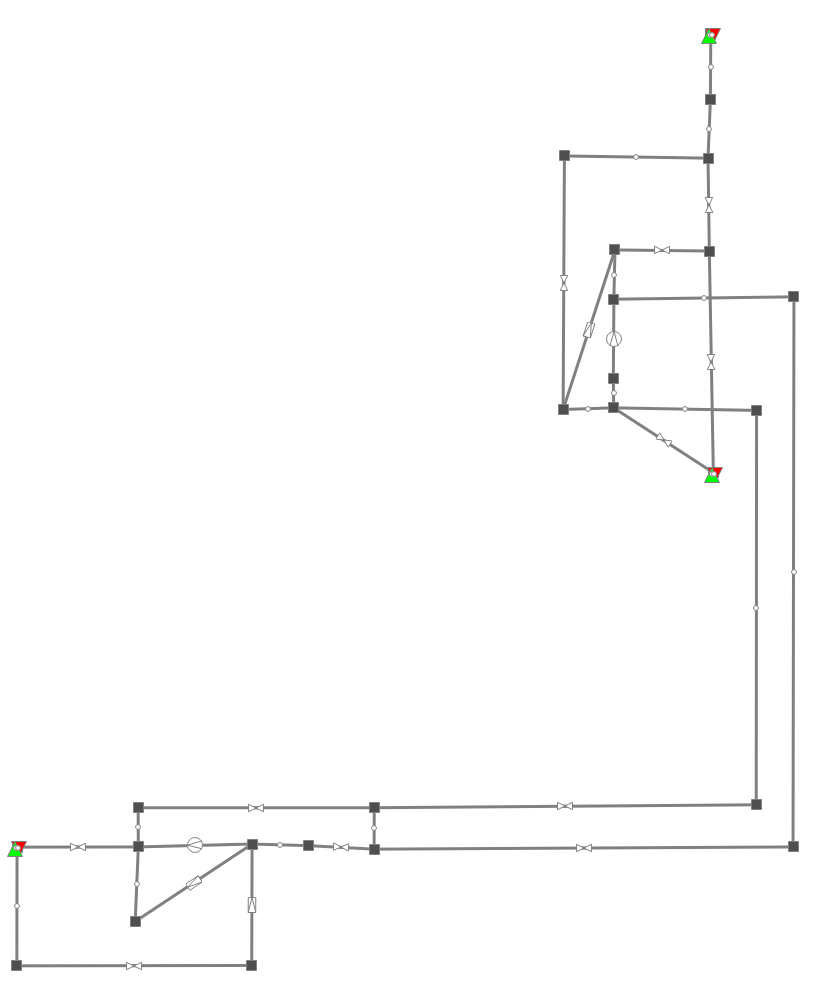}
    \caption{Topology of Station D.}
    \label{fig:porz_topology}
\end{figure}

All instances, both in training and test, contain 12 time steps (excluding the initial state) with 30 minutes between each step. Additionally, we focus on Station D from \cite{hennings2020controlling}, and present only results for this station. The statistics for Station D can be seen in Table \ref{tab:stationStatistics}, and its topology in Figure \ref{fig:porz_topology}. Station D can be thought of as a T intersection, and is of average complexity compared to the stations presented in \cite{hennings2020controlling}. The station contains 6 boundary nodes, but they are paired, such that for each pair only one can be active, i.e., have non-zero flow. Due to this, our sampling method in subsection \ref{subsection:Boundary Forecast Generation} exists in 3-dimensions and is uniform $\forall t \in \setTimestepsNoZero$. 

\begin{table}[ht]
    \centering
    \begin{tabular}{*{8}{c}}
      Name & $|\setVertices|$ & $|\setArcs|$ & $\frac{\sum\limits_{a\in\setPipes}\pipeLength{a}}{|\setPipes|}$ &
      $|\setCompressorConfigurations{a}| \enskip \forall a\in\setCompressorStations$ &
      $|\setNSModes|$ & $|\setBoundaryNodes|$ & $|\setValves|$ \\
      \midrule
      D              &  31 &  37 & 0.404 km &           2, 6 &   56  & 3x2 & 11 \\
    \end{tabular}
    \caption{Overview of different properties of station D.}
    \label{tab:stationStatistics}
\end{table}

\subsection{Exact Network Designs}
As a large portion portion of our input data into both \generator and \discriminator is time-expanded data, we originally believed that the ideal design would be a series of LSTMs \cite{hochreiter1997long}. Preliminary results however showed that convolutional neural networks (CNNs) were more effective for our problem, in particular when using Inception Blocks, see \cite{szegedy2017inception}. 

The exact block design used in \neuralnetwork can be seen in Figure \ref{fig:block_design}, and the general layout in Figure \ref{fig:forward_pass}. For the complete network design we refer readers to Figure \ref{fig:graph_neural_network} and Table \ref{tab:nn_architecture} in the Appendix.

%% file: sections/TrainingAlgo.tex
\begin{algorithm}[H]
\SetAlgoLined
\SetKwInOut{Input}{Input}
\Input{Neural network \neuralnetwork, prelabelled\_data }
\KwResult{Trained neural network \neuralnetwork}
set\_trainable(\discriminator)\; 
set\_untrainable(\generator)\;
Discriminator Pretraining(\discriminator, prelabelled\_data) \footnote{See Algorithm \ref{alg:pretrain-discriminator}} \;
softmax\_temperature = 0\;
data = []\;
\For{$i = 0;~i < \text{num\_epochs}$}{
set\_trainable(\generator)\; 
set\_untrainable(\discriminator)\;
\For{$i = 0;~i < \text{num\_generator\_epochs}$}{
softmax\_temperature += 1\;
set(\generator, softmax\_temperature)\;
loss = Generator Training(\neuralnetwork) \footnote{See Algorithm \ref{alg:train-generator}} \;
\If{$\text{loss} \leq \text{stopping\_loss\_generator}$}{
break\;
}
}
set\_trainable(\discriminator)\; 
set\_untrainable(\generator)\;
data = Prepare Discriminator Training Data(\neuralnetwork, data) \footnote{See Algorithm \ref{alg:data_prepare}} \;
mixed\_data = MixData(data, prelabelled\_data, num\_prelabelled)\;
training\_data, test\_data = split\_data(mixed\_data, ratio\_test)\;
optimizer = Adam(learning\_rate, weight\_decay) \footnote{See \cite{kingma2014adam}  \href{pytorch.org/docs/stable/optim.html?highlight=adam\#torch.optim.Adam}{pytorch.org/docs/stable/optim.html?highlight=adam\#torch.optim.Adam}.} \;
lr\_scheduler=ReduceLROnPlateau \footnote{See \href{pytorch.org/docs/stable/optim.html\#torch.optim.lr\_scheduler.ReduceLROnPlateau} {pytorch.org/docs/stable/optim.html\#torch.optim.lr\_scheduler.ReduceLROnPlateau}.}(patience, factor)\;
dataloader = DataLoader(training\_data, batch\_size, shuffle=True)\;
\For{$i = 0;~i < \text{num\_discriminator\_epochs}$}{
Discriminator Training Loop(\discriminator, dataloader, optimizer) \footnote{See Algorithm \ref{alg:training-loop-discriminator}} \;
lr\_scheduler.step()\;
test\_loss = compute\_L1Loss(\discriminator, test\_data)\;
\If{$\text{test\_loss} \leq \text{stopping\_loss\_discriminator}$}{
break\;
}
}
}
\Return \neuralnetwork
\caption{Neural Network Training}
\label{alg:train}
\end{algorithm}

%% file: sections/Results.tex
We partition our results into three subsections. The first focuses on the training results of \neuralnetwork, the second on our data generation methods, while the third is concerned with our results on the 15 weeks of real-world transient gas data. Note that when training we scaled \objfixedmip values by 500 to reduce the magnitude of the losses. For visualisation purposes of comparing the performance of \neuralnetwork and our data generation methods, we re-scaled all results.

\subsection{Training Results}
\begin{figure}[ht]
    \centering
    \includegraphics[scale=0.75]{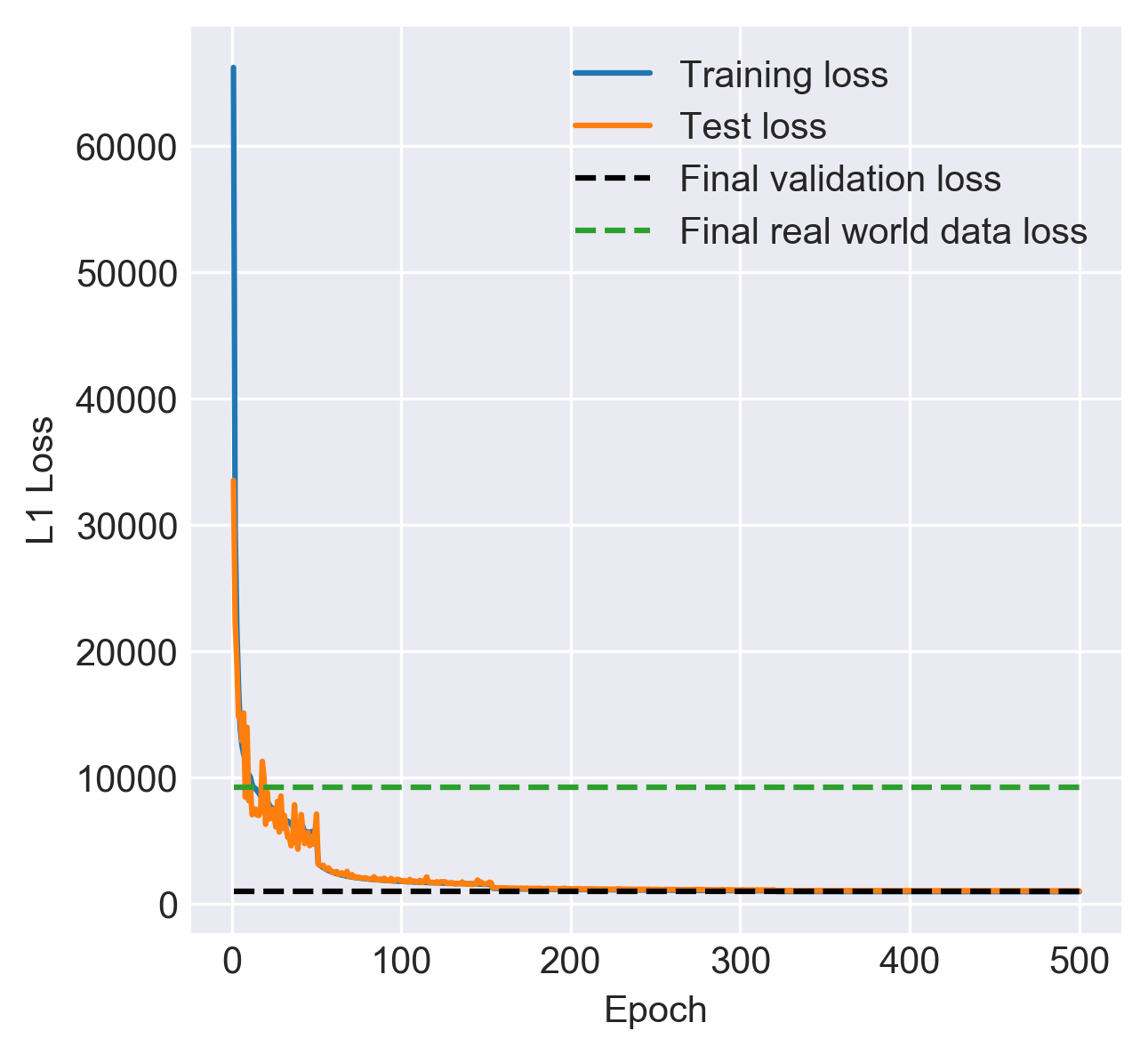}
    \caption{The loss per epoch of \discriminator during the initial training of Algorithm \ref{alg:pretrain-discriminator}. The dashed lines show the performance of \discriminator after \neuralnetwork has been completely trained.}
    \label{fig:training_discriminator}
\end{figure}

Figure \ref{fig:training_discriminator} shows the training loss throughout the initial offline training. We see that \discriminator learns how to accurately predict \objfixedmip as the loss decreases. This is a required result, as without a trained discriminator we cannot expect to train a generator. Both the training and test loss converge to approximately 1000, which is excellent considering the generated \objfixedmip range well into the millions. As visible by both the test loss and final validation loss, we see \discriminator generalises to \fixedmip instances of our validation set that it has not seen. This generalisation ability doesn't translate perfectly to real-world data however. This is due to the underlying distribution of real-world data and our generated data being substantially different. Despite this we believe that an L1 loss, in this case simply the average distance between \predobjfixedmip and \objfixedmip, of 10000 is still very good. We discuss the issues of different distributions in subsection \ref{sec:data_generation_results}.

\begin{figure}[h]
\centering
\includegraphics[scale=0.65]{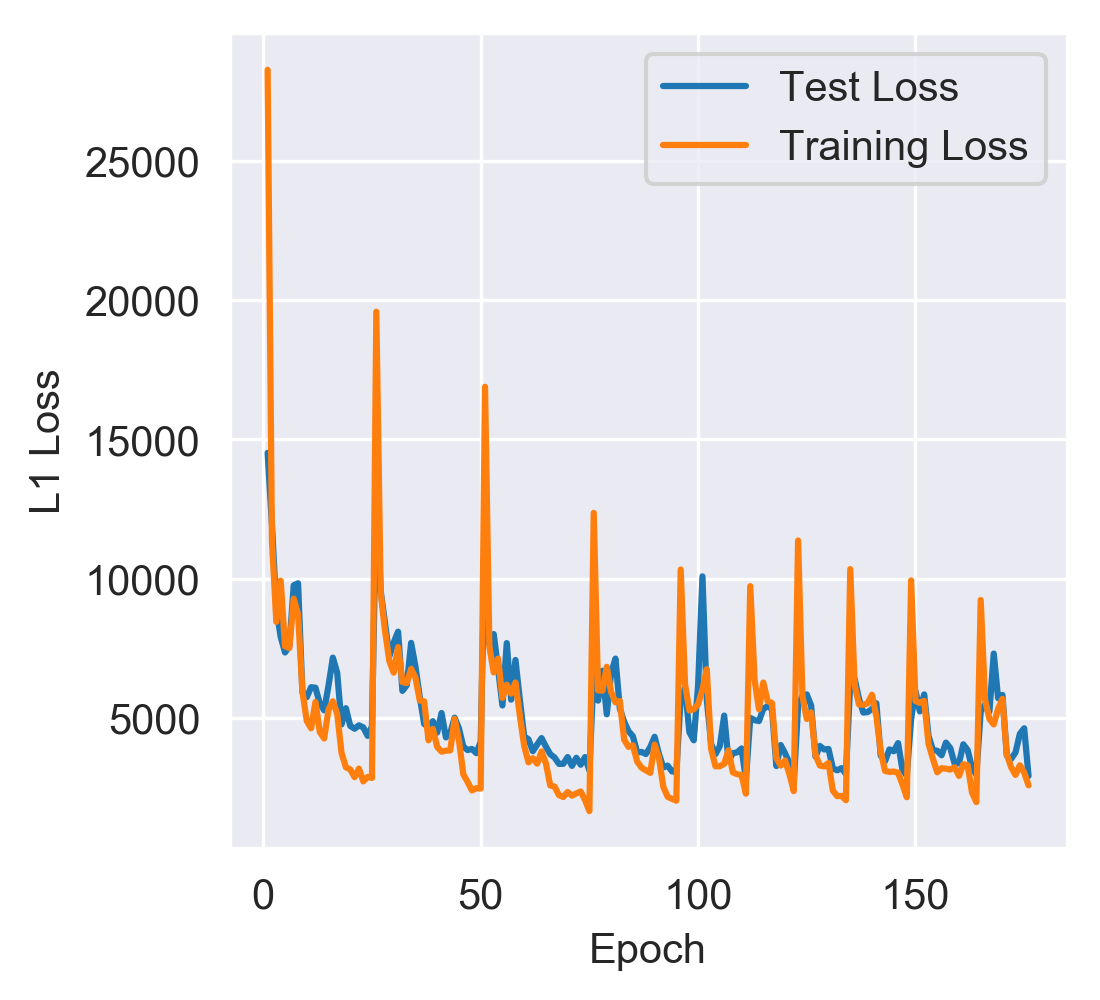}
\caption{The loss per epoch of \discriminator as it is trained using Algorithm \ref{alg:train}}
\label{fig:retraining_discriminator}
\end{figure}

\begin{figure}[h]
    \centering
    \includegraphics[scale=0.60]{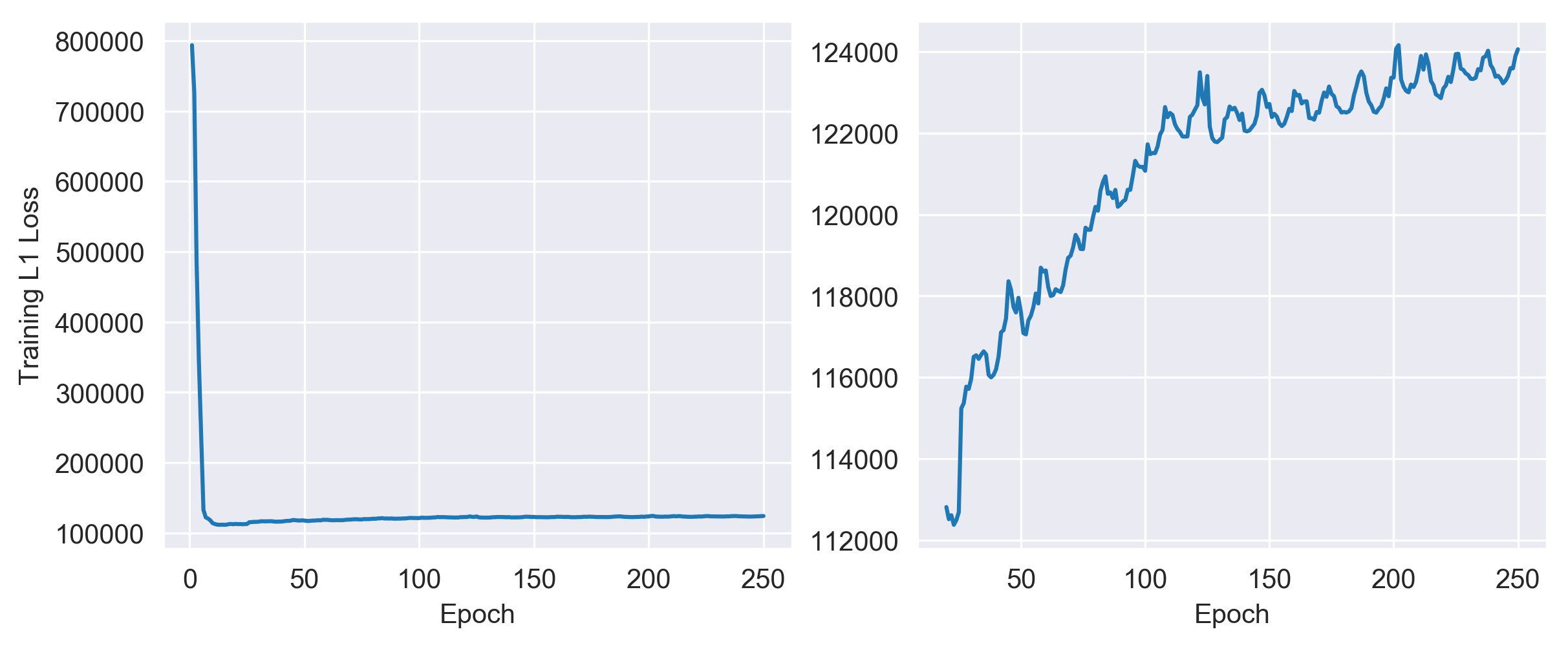}
    \caption{(Left) The loss per epoch of \generator as it is trained using Algorithm \ref{alg:train}. On the left the loss over all epochs is shown. (Right) A magnified view of the loss starting from epoch 20.}
    \label{fig:full_training_generator}
\end{figure}

The loss during training using Algorithm~\ref{alg:train} for \discriminator is shown in Figure \ref{fig:retraining_discriminator}, and for \generator in Figure \ref{fig:full_training_generator}. The cyclical nature of the \discriminator loss is caused by the re-training of \generator, which learns how to induce sub-optimal predictions from the then static \discriminator. These sub-optimal predictions are quickly re-learned, but highlight that learning how to perfectly predict \objfixedapproxmip over all possibilities, potentially due to the rounded nature of \predbinvars, is unlikely without some error. Figure \ref{fig:full_training_generator} (left) shows the loss over time of \generator as it is trained, with Figure \ref{fig:full_training_generator} (right) displaying magnified losses for the final epochs. We observe that \generator quickly learns important \binvars decision values. We hypothesise that this quick descent is helped by \predbinvars that are unlikely given our generation method in Algorithm \ref{alg:op_mode_generator}. The loss increases following this initial decrease in the case of \generator, showing the ability of \discriminator to further improve. It should also be noted that significant step-like decreases in loss are absent in both (left) and (right) of Figure \ref{fig:full_training_generator}. Such steps would indicate \generator discovering new important \binvars values (operation modes). The diversity of produced operation modes however, see Figure \ref{fig:generated_opmodes}, implies that early in training a complete spanning set of operation modes is derived, and the usage of their ratios is then learned and improved.

\subsection{Data Generation Results}
\label{sec:data_generation_results}

\begin{figure}[h]
    \centering
    \includegraphics[scale=0.6]{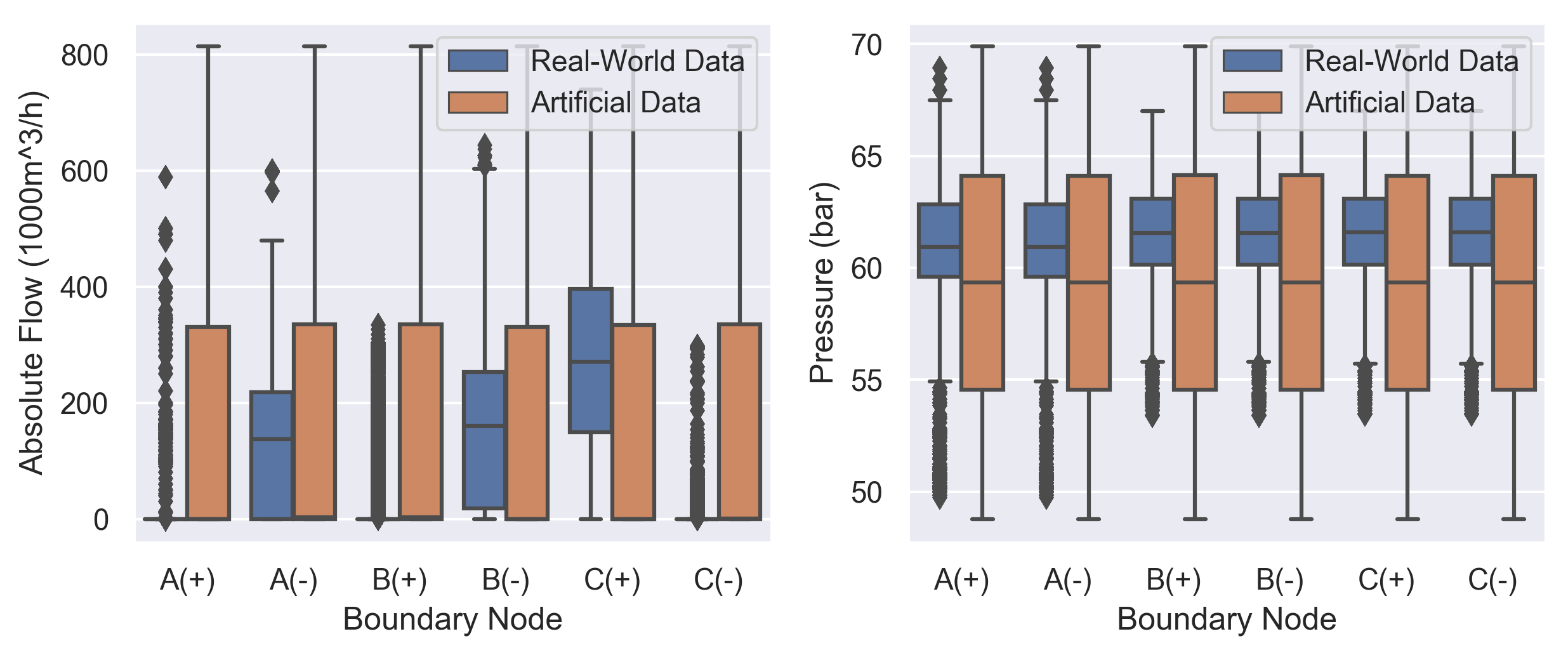}
    \caption{Comparison of generated flow (Left) / pressure (Right) value distributions per node vs. the distribution seen in real-world data.}
    \label{fig:generated_boundary_values}
\end{figure}

As an interlude between results from \neuralnetwork, we outline the performance of our synthetic gas network data generation methods. Figure \ref{fig:generated_boundary_values} (left) shows how our generated flow prognosis compares to that of historic real-world data. We see that Nodes A, B, and C are not technically entry or exits, but over historical data are dominated by a single orientation for each node. Specifically, Node C is the general entry, and Nodes A / B are the exits. In addition to the general orientation, we see that each node has significantly different ranges and distributions. These observations highlight the simplicity of our data generation methods, as we see near identical distributions for all nodes over the artificial data. We believe this calls for further research in prognosis generation methods. Figure \ref{fig:generated_boundary_values} (right) shows our pressure prognosis compared to that of historic values. Unlike historic flow values, we observe little difference between historic pressure values of different nodes. This is supported by the optimal choices \optbinvars over the historic data, see Figure \ref{fig:generated_opmodes}, as in a large amount of cases compression is not needed and the network station is in bypass. Note that each corresponding entry (+) and exit (-) have identical pressure distributions due to the way they are constructed.

\begin{figure}[h]
    \centering
    \includegraphics[scale=0.6]{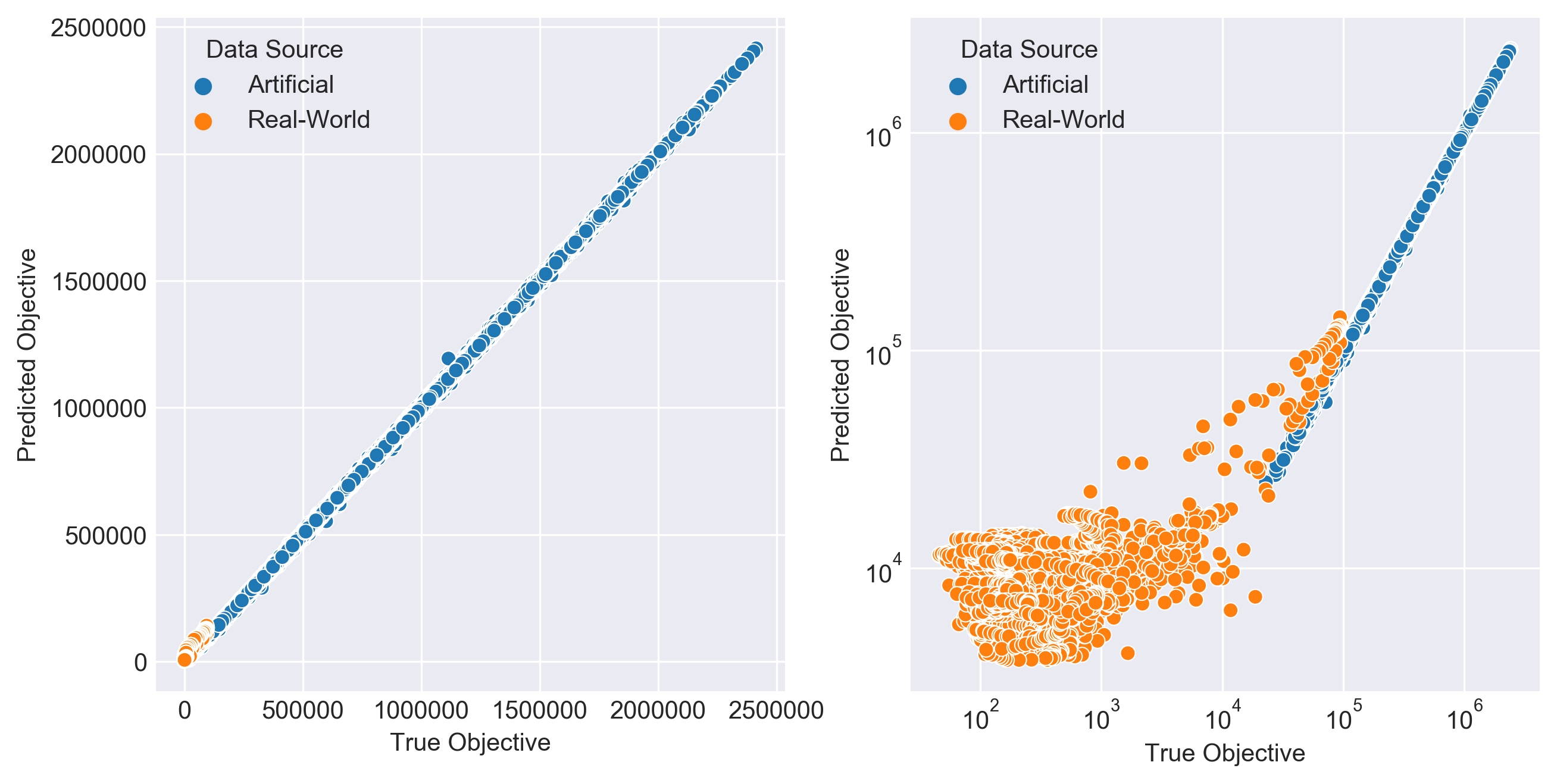}
    \caption{\predobjfixedapproxmip for the validation set, and \predobjmip for real-world data, compared to \objfixedapproxmip and \objmip respectively. Linear scale (Left) and log-scale (Right). }
    \label{fig:pred_obj_trained}
\end{figure}

A further comparison of how our generated data compares to historic data can be seen in Figure \ref{fig:pred_obj_trained}. Here one can see the distribution of \predobjfixedapproxmip and \objfixedapproxmip for the generated validation set, and \predobjmip and \objmip for the real-world data. As expected, the distributions are different depending on whether the data is artificial or not. Our data generation was intended to be simplistic, and as independent as possible from the historic data. As such, the average scenario has optimal solution larger than that of any real-world data point. The performance of \discriminator is again clearly visible here, with \predobjfixedapproxmip and \objfixedapproxmip being near identical over the artificial data, keeping in mind that these data points were never used in training. We see that this ability to generalise is relatively much worse on real-world data, mainly due to the the lower values of \objmip over this data. Figure \ref{fig:pred_obj_trained} (right) shows the results with log-scale axes to better highlight this disparity. It should be noted that the real-world instances with larger \objmip are predicted quite well, and all real-world instances have an L1 distance between \predobjmip and \objmip that is small in terms of absolute differences. 

\subsection{Real-World Results}

\begin{figure}
    \centering
    \includegraphics[scale=0.6]{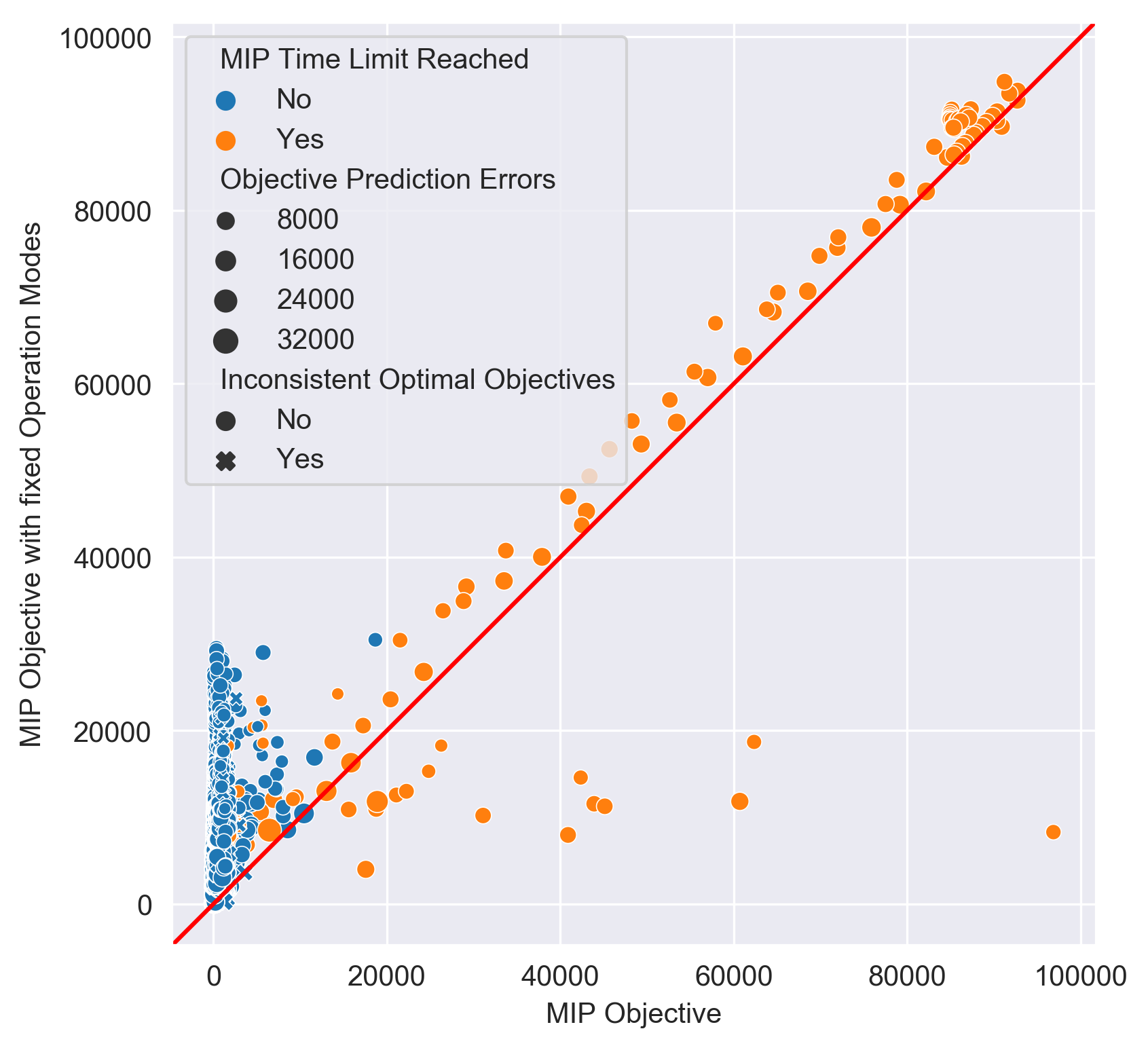}
    \caption{A comparison of \objfixedapproxmip and \objmip for all real-world data instances.}
    \label{fig:fixed_obj_vs_true_obj}
\end{figure}

We now present results of our fully trained \neuralnetwork applied to the 15 weeks of real-world data. Note that we had to remove 651 instances from our 9291 instances, as the warm-start resulted in an optimal solution value further away than the optimality tolerances we set. These instances have been kept in the graphics, but are marked and conclusions will not be drawn from them. We believe the problems with reproducibility are caused by the numeric difficulties in managing the pipe equality constraints. 

Figure \ref{fig:fixed_obj_vs_true_obj} shows the comparison of \objfixedapproxmip and \objmip. In a similar manner to \discriminator, we see that \generator struggles with instances where \objmip is small. This is visible in the bottom left, where we see \objfixedapproxmip values much larger than \objmip for like \instance. This comes as little surprise given the struggle of \discriminator with small \objmip values. Drawing conclusions becomes more complicated for instances with larger \objmip values, because the majority hit the time limit. We can clearly see however, the value of our primal heuristic. There are many cases, those below the line \objfixedapproxmip = \objmip, where our primal heuristic retrieves a better solution than the MILP solver does in one hour. Additionally, we see that no unsolved point above the line is very far from the line, showing that our primal heuristic produced a comparable, sometimes equivalent solution in a much shorter time frame. For a comparison of solve-times, see Table \ref{tab:solve_times}. 

\begin{table}[!ht]
\centering
\begin{small}
\begin{tabular}{lrrrrr}
\toprule
{} &     Mean &  Median &      STD &    Min &       Max \\
\midrule
\neuralnetwork Inference Time (s)               &    0.009 &   0.008 &    0.001 &  0.008 &     0.017 \\
Warmstarted \mip Time (s)                        &  100.830 &   9.380 &  421.084 &  0.130 &  3600.770 \\
\mip Time (s)                                    &  147.893 &  24.380 &  463.279 &  3.600 &  3601.280 \\
\fixedapproxmip + Warmstarted \mip Time (s) &  103.329 &  12.130 &  424.543 &  0.190 &  3726.110 \\
\fixedapproxmip Time (s)                        &    2.499 &   1.380 &   12.714 &  0.060 &   889.380 \\
\bottomrule
\end{tabular}
\end{small}
\caption{Solve time statistics for different solving strategies.}
\label{tab:solve_times}
\end{table}

\begin{figure}[!ht]
    \centering
    \includegraphics[scale=0.6]{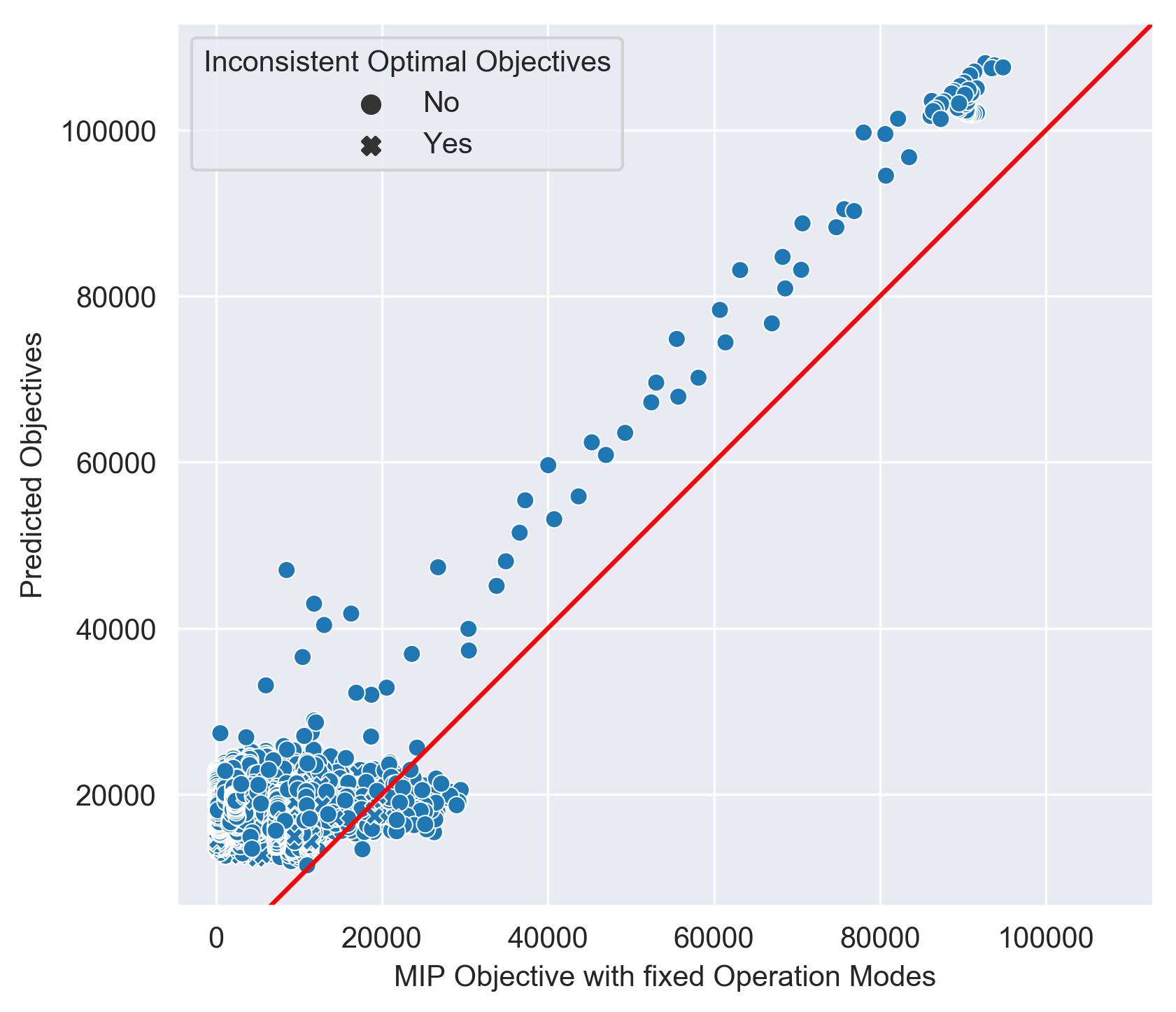}
    \caption{A comparison of \predobjfixedapproxmip and \objfixedapproxmip for all real-world data instances.}
    \label{fig:fixed_obj_vs_pred_obj}
\end{figure}

\newpage

Figure \ref{fig:fixed_obj_vs_pred_obj} shows the performance of the predictions \predobjfixedapproxmip compared to \objfixedapproxmip. Interestingly, \discriminator generally predicts \predobjfixedapproxmip values slightly larger than \objfixedapproxmip. We expect this for the smaller valued instances, as we know that \discriminator struggles with \objfixedapproxmip instances near 0, but the trend is evident for larger valued instance too. The closeness of the data points to the line \predobjfixedapproxmip = \objfixedapproxmip show that \discriminator can adequately predict \predbinvars solutions from \generator despite the change in data sets. Figure \ref{fig:fixed_obj_vs_true_obj} showed that \generator successfully generalised to a new data set, albeit with difficulties around instances with \objmip valued near 0.  From Figures \ref{fig:fixed_obj_vs_true_obj} and \ref{fig:fixed_obj_vs_pred_obj}, we can see that the entire \neuralnetwork generalises to unseen real-world instances, despite some generalisation loss. 

\begin{figure}
    \centering
    \includegraphics[scale=0.6]{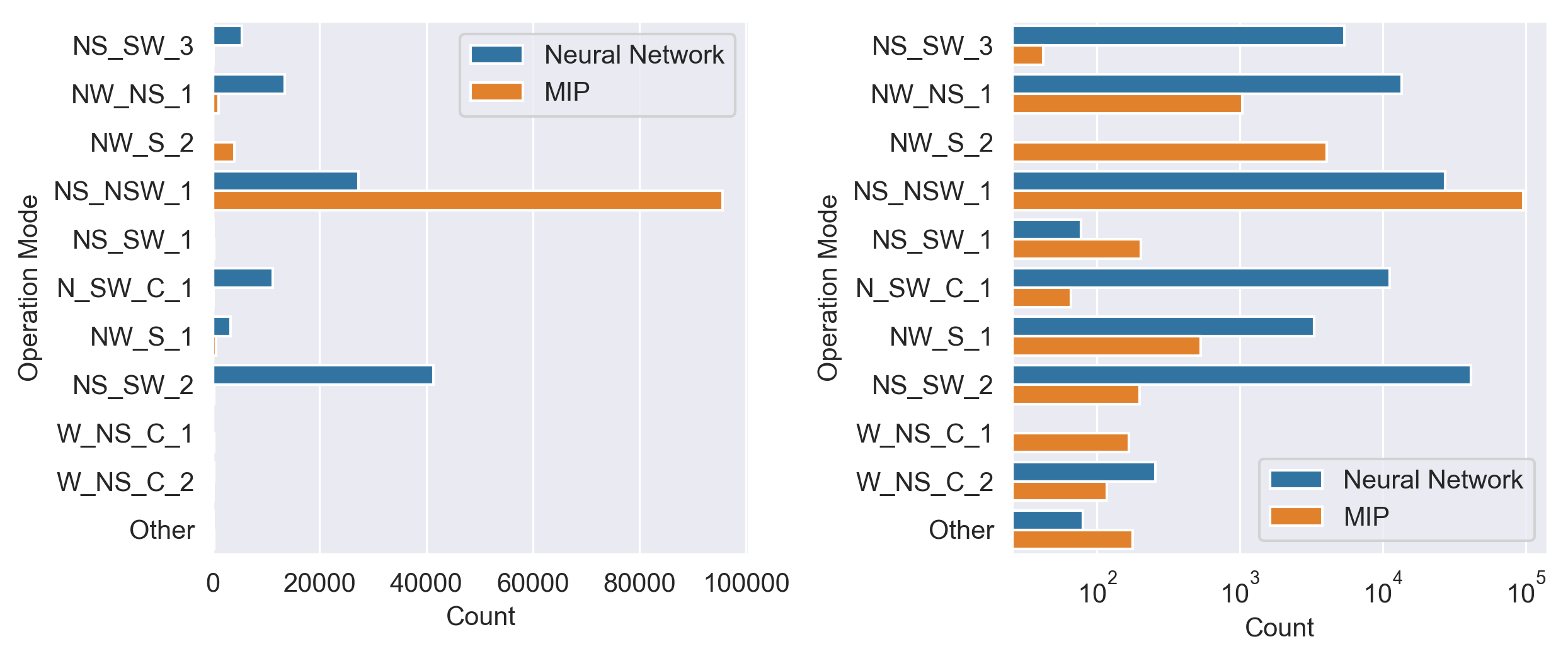}
    \caption{Frequency of operation mode choice by \generator compared to MILP solver for all real-world instances. (Left) Linear scale, and (Right) log scale.}
    \label{fig:generated_opmodes}
\end{figure}

\input{sections/Pictures/operation_modes_matrix}

We now compare the operation modes \predbinvars, which are generated by \generator, and the \optbinvars, which are produced by our MILP solver. To do so we use the following naming convention: We name the three pairs of boundary nodes N (north), S (south), and W (west). Using W\_NS\_C\_2 as an example, we know that flow comes from W, and goes to N and S. The C in the name stands for active compression, and the final index is to differentiate between duplicate names. As seen in Figure \ref{fig:generated_opmodes}, which plots the frequency of specific \binvars if they occurred more than 50 times, a single choice dominates \optbinvars. This is interesting, because we expected there to be a-lot of symmetry between \binvars, with the MILP solver selecting symmetric solutions with equal probability. For instance, take W\_NS\_C\_1 and take W\_NS\_C\_2. \neuralnetwork only ever predicts W\_NS\_C\_2, however with half the frequency the MILP solver selects each of them. This indicates that from the MILP's point of view they are symmetric, and either can be chosen, while \neuralnetwork has recognised this and converged to a single choice. We can support this by analysing the data, where the difference in W\_NS\_C\_1 and W\_NS\_C\_2 is which compressor machine is used, with both machines being identical. This duplicate choice apparently does not exist in bypass modes however, where the uniqueness of \binvars, determined by valve states, results in different \objfixedmip values. It is observable then that for the majority of instances NS\_NSW\_1 is the optimal choice, and that \neuralnetwork has failed to identify its central importance. We believe this is due to the training method, where over generalisation to a single choice is strongly punished. For a comprehensive overview of the selection of operation modes and the correlation between \predbinvars and \optbinvars, we refer interested readers to Table \ref{tab:operation_mode_correlation}.

\begin{figure}[ht]
    \centering
    \includegraphics[scale=0.6]{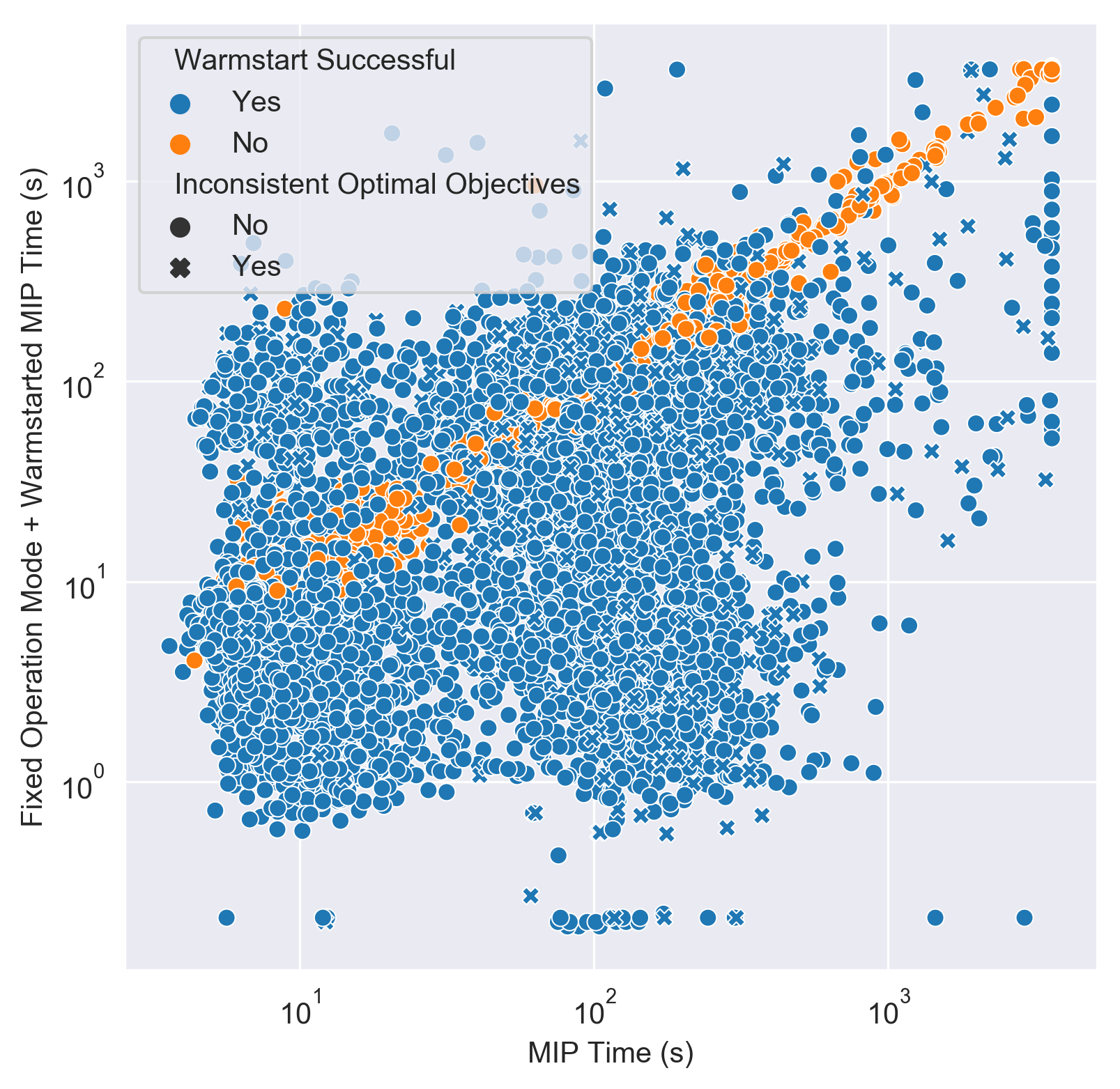}
    \caption{The combined running time of solving \fixedapproxmip, and solving a warm-started \mip, compared to solving \mip directly.}
    \label{fig:warm_start_plus_lp}
\end{figure}

As discussed above, \neuralnetwork cannot reliably produce \optbinvars. Nevertheless, it produces near-optimal \predbinvars suggestions, which are still useful in a warm-start context, see Algorithm \ref{alg:warm_start}. The results of our warm-start algorithm are displayed in Figure \ref{fig:warm_start_plus_lp}. Our warm-start suggestion was successful 72\% of the time, and the algorithm resulted in an average speed up of 60.5\%. We use the shifted geometric mean with a shift of 1 for this measurement to avoid distortion by relative variations of the smaller valued instances. Especially surprising is that some instances that were previously unsolvable within the time-limit were easily solvable given the warm-start suggestion. In addition, many of the solvable but complicated instances are also solved near instantly with the warm-start suggestion. As such, we have created an effective primal heuristic that is both quick to run and beneficial in the context of locating a globally optimal solution. 

%% file: sections/Pictures/operation_modes_matrix.tex
\begin{table}
    \centering
    \begin{tiny}
    \bgroup
    \def\arraystretch{1.5}
    \rowcolors{2}{gray!25}{white}
    \begin{tabular}{lrrrrrrrrrrr}
    {} & \rotatebox{70}{NW\_NS\_1} & \rotatebox{70}{NS\_SW\_2} & \rotatebox{70}{N\_SW\_C\_1} &\rotatebox{70}{NS\_NSW\_1}  &\rotatebox{70}{W\_NS\_C\_1}  & \rotatebox{70}{NS\_SW\_1} & \rotatebox{70}{NW\_S\_2} &\rotatebox{70}{NS\_SW\_3}  &\rotatebox{70}{W\_NS\_C\_2}  & \rotatebox{70}{NW\_S\_1} & \rotatebox{70}{Other} \\
    \midrule
    NW\_NS\_1  &     884 &      22 &        0 &     9529 &       31 &      37 &   2436 &       4 &       24 &    397 &    82 \\
    NS\_SW\_2  &      48 &     102 &        1 &    40298 &        0 &     114 &    630 &      24 &        0 &     51 &    13 \\
    N\_SW\_C\_1 &       0 &      27 &       65 &    11008 &        0 &       4 &      0 &       2 &        0 &      0 &    55 \\
    NS\_NSW\_1 &      41 &      29 &        0 &    26509 &        0 &      28 &    557 &       9 &        0 &     49 &    15 \\
    W\_NS\_C\_1 &       0 &       0 &        0 &        0 &        0 &       0 &      0 &       0 &        0 &      0 &     0 \\
    NS\_SW\_1  &       0 &       0 &        0 &       76 &        0 &       1 &      0 &       0 &        0 &      0 &     0 \\
    NW\_S\_2   &       4 &       0 &        0 &        0 &        0 &       0 &      2 &       0 &        0 &      1 &     1 \\
    NS\_SW\_3  &       6 &       7 &        0 &     5220 &        0 &       7 &    108 &       1 &        0 &      4 &     5 \\
    W\_NS\_C\_2 &      28 &       0 &        0 &        0 &      136 &       0 &      0 &       0 &       93 &      0 &     0 \\
    NW\_S\_1   &      30 &      11 &        0 &     2880 &        0 &      12 &    315 &       2 &        0 &     30 &     6 \\
    Other    &       0 &       1 &        0 &       78 &        0 &       0 &      0 &       0 &        0 &      0 &     1 \\
    \bottomrule
    \end{tabular}
    \egroup
    \end{tiny}
\caption{Operation Mode Correlation Matrix between \predbinvars and \optbinvars.}
\label{tab:operation_mode_correlation}
\end{table}

%% file: sections/Conclusion.tex
In this paper, we presented a dual neural network design for generating decisions in a MILP. This design is trained without ever solving the MILP with unfixed decision variables. The neural network is both used as a primal heuristic and used to warm-start the MILP solver for the original problem. We proved the usefulness of our design on the transient gas transportation problem. While doing so we created methods for generating synthetic transient gas data for training purposes, reserving an unseen 9291 real-world instances for validation purposes. Despite some generalisation loss, our trained neural network results in a primal heuristic that takes on average 2.5s to run, and results in a 60.5\% decrease in global optimal solution time when used in a warm-start context.

While our approach is an important step forward in neural network design and ML's application to gas transport, we believe that there exists four primary directions for future research. The first of which is to convert our approach into more traditional reinforcement learning, and then utilise policy gradient approaches, see \cite{thomas2017policy}. The major hurdle to this approach is that much of the computation would be shifted online, requiring many more calls to solve the induced MILPs. This could be offset however, by using our technique to initialise the weights for such an approach, thereby avoiding early stage training difficulties with policy gradient approaches. The second is focused on the recent improvements in Graph Neural Networks, see \cite{gasse2019exact}. Their ability to generalise to different input sizes would permit the creation of a single NN over multiple network stations or gas network topologies. Thirdly, there exists a large gap in the literature w.r.t data generation for transient gas networks. Improved methods are needed, which are scalable and result in real-world like data. Finally, although we focused on the transient gas transportation problem, our approach can be generalised to arbitrary problem classes. 

%% file: sections/Funding.tex
The work for this article has been conducted in the Research Campus MODAL funded by the German Federal Ministry of Education and Research (BMBF) (fund numbers 05M14ZAM, 05M20ZBM), and was supported by the German Federal Ministry of Economic Affairs and Energy (BMWi) through the project UNSEEN (fund no 03EI1004D).

%% file: sections/appendix/appendix.tex
\section{Appendix} \label{sec:appendix}
\input{sections/appendix/DataPrepareAlgo}
\input{sections/appendix/PretrainDiscriminator}
\input{sections/appendix/DiscriminatorTraining}
\input{sections/appendix/GeneratorTraining}
\input{sections/appendix/TableInitialTraining}
\input{sections/appendix/MIPParams}
\begin{figure}
    \centering
    \includegraphics[scale=0.55]{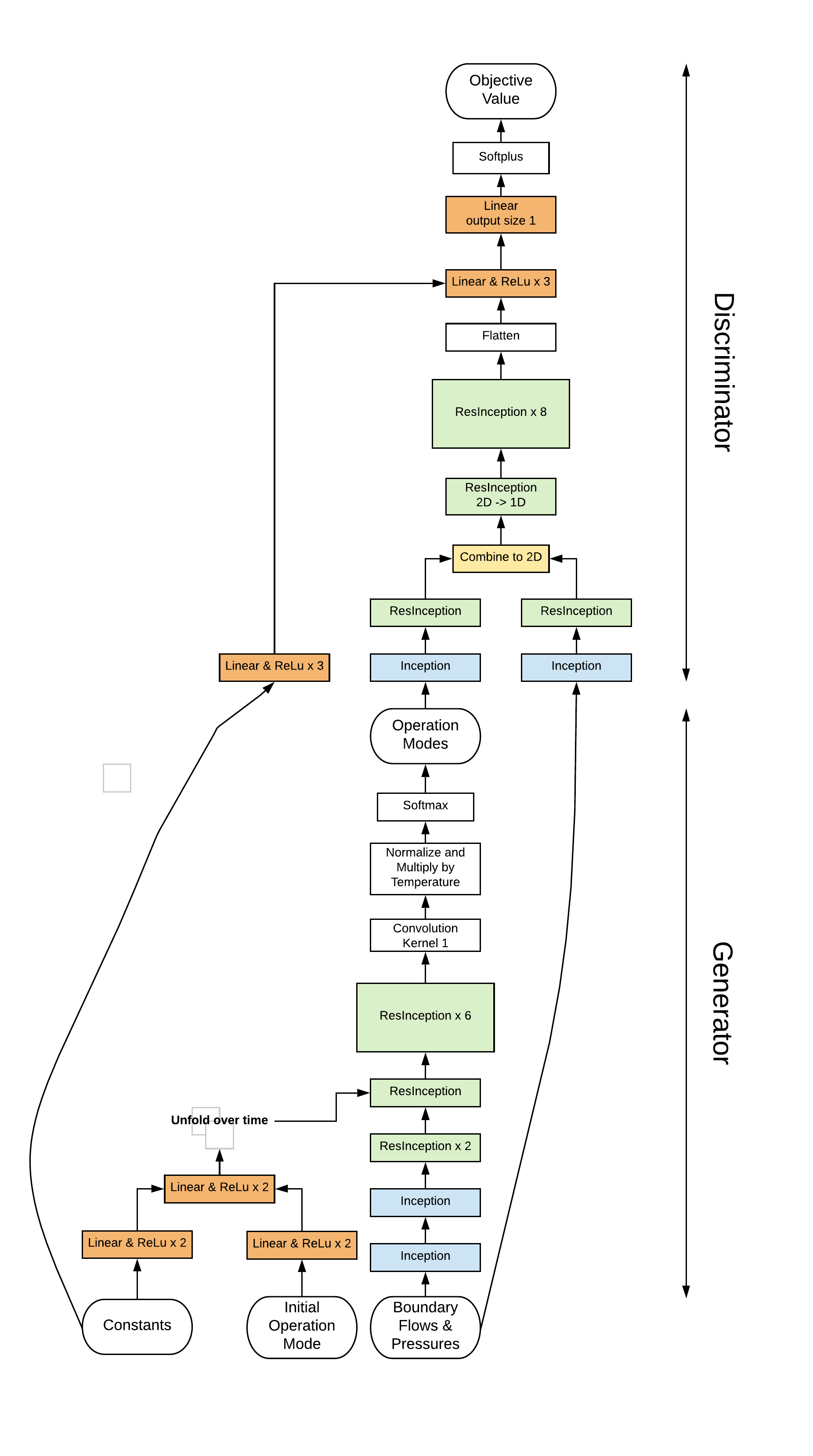}
    \caption{Neural Network Architecture}
    \label{fig:graph_neural_network}
\end{figure}
\input{sections/appendix/parameterTable}

%% file: sections/appendix/DataPrepareAlgo.tex
\begin{algorithm}[H]
\SetAlgoLined
\SetKwInOut{Input}{Input}
\Input{Neural network \neuralnetwork , labelled\_data}
\KwResult{Data for training \discriminator}
new\_labelled\_data = [] \;
\For{$i = 0;\ i < \text{num\_data\_new}$}{
initial\_state = Uniformly select from generated offline data \;
flow\_forecast, pressure\_forecast = Boundary Prognosis Generator()\footnote{See Algorithm \ref{alg:boundary_data_generator}}\;
\instance = (flow\_forecast, pressure\_forecast, initial\_state)\;
\predbinvars, \predobjfixedapproxmip = \neuralnetwork($\pi$)\;
\objfixedapproxmip = solve \fixedapproxmip \;
new\_labelled\_data.append(\instance, \predbinvars, \objfixedapproxmip) \;
}

\Return concatenate(labelled\_data[-num\_data\_old:], new\_labelled\_data)
\caption{Prepare Discriminator Training Data}
\label{alg:data_prepare}
\end{algorithm}

%% file: sections/appendix/PretrainDiscriminator.tex
\begin{algorithm}[h]
\SetAlgoLined
\SetKwInOut{Input}{Input}
\Input{Discriminator \discriminator, data}
optimizer = Adam(learning\_rate, weight\_decay)\;
dataloader = DataLoader(data, batch\_size, shuffle=True)\;
lr\_scheduler=ReduceLROnPlateau()\;
\For{$i = 0;~i < \text{num\_epochs}$}{
Discriminator Training Loop(\discriminator, dataloader, optimizer) \footnote{See Algorithm \ref{alg:training-loop-discriminator}}\; 
lr\_scheduler.step()\;
}
\caption{Discriminator Pretraining}
\label{alg:pretrain-discriminator}
\end{algorithm}

%% file: sections/appendix/DiscriminatorTraining.tex
\begin{algorithm}
\SetAlgoLined
\SetKwInOut{Input}{Input}
\Input{Discriminator \discriminator, dataloader, optimizer}
\For{batch in dataloader}{
optimizer.zero\_grad()\;
\predobjfixedmip = \discriminator(batch)\;
loss = L1Loss(\predobjfixedmip, \objfixedmip)\;
loss.backward()\;
optimizer.step()\;
}
\caption{Discriminator Training Loop}
\label{alg:training-loop-discriminator}
\end{algorithm}

%% file: sections/appendix/GeneratorTraining.tex
\begin{algorithm}[H]
\SetAlgoLined
\SetKwInOut{Input}{Input}
\Input{Neural network \neuralnetwork}
\KwResult{Average loss in training}
optimizer = Adam()\;
lr\_scheduler=cyclicLR\footnote{Introduced by Smith in \cite{smith2017cyclical}, see also \newline \href{pytorch.org/docs/stable/optim.html\#torch.optim.lr\_scheduler.CyclicLR} {pytorch.org/docs/stable/optim.html\#torch.optim.lr\_scheduler.CyclicLR}.}(max\_lr, base\_lr, step\_size\_up)\;
data = [] \;
\For{$i = 0;\ i < \text{num\_scenarios}$}{
initial\_state = Uniformly select from generated offline data \;
flow\_forecast, pressure\_forecast = Boundary Prognosis Generator()\footnote{See Algorithm \ref{alg:boundary_data_generator}}\;
\instance = (flow\_forecast, pressure\_forecast, initial\_state)\;
data.append(\instance) \;
}
dataloader = DataLoader(data, batch\_size, shuffle=True)\;
\For{batch in dataloader}{
optimizer.zero\_grad()\;
losses = []\;
\predbinvars, \predobjfixedapproxmip = \neuralnetwork(batch)\;
loss = L1Loss(\predobjfixedapproxmip, 0)\;
loss.backward()\;
optimizer.step()\;
lr\_scheduler.step()\;
losses.append(loss)\;
}
\KwResult{mean(losses)}
\caption{Generator Training}
\label{alg:train-generator}
\end{algorithm}

%% file: sections/appendix/TableInitialTraining.tex
\begin{table}[!ht]
    \centering
    \begin{threeparttable}[t]
    \begin{tabular}{lll}
        \toprule
        Parameter               & Method                                            & Value \\
        \midrule
        batch\_size                 &Algorithm~\ref{alg:pretrain-discriminator}         &   2048 \\
        num\_epochs                 &Algorithm~\ref{alg:pretrain-discriminator}         &    500 \\
        learning\_rate              &Algorithm~\ref{alg:pretrain-discriminator} / Adam  &  0.005 \\
        weight\_decay               &Algorithm~\ref{alg:pretrain-discriminator} / Adam  &  5e-06 \\
        
        batch\_size                 &Algorithm~\ref{alg:train-generator}                &   2048 \\
        max\_lr                     &Algorithm~\ref{alg:train-generator} / CyclicLR     &  0.0005 \\
        base\_lr                    &Algorithm~\ref{alg:train-generator} / CyclicLR     &  5e-06 \\
        step\_size\_up              &Algorithm~\ref{alg:train-generator} / CyclicLR     & 10000\\
        num\_scenarios              &Algorithm~\ref{alg:train-generator}                & 3200000 \\
        num\_data\_new   &Algorithm~\ref{alg:data_prepare}                   & 2048 \\
        num\_data\_old        &Algorithm~\ref{alg:data_prepare}                   & 8192 \\
        num\_epochs                 &Algorithm~\ref{alg:train}                          &10 \\
        num\_generator\_epochs      &Algorithm~\ref{alg:train}                          &25 \\
        num\_discriminator\_epochs  &Algorithm~\ref{alg:train}                          &25 \\
        stopping\_loss\_discriminator&Algorithm~\ref{alg:train}                         &3 * 1022.5\tnote{1}\\
        stopping\_loss\_generator   &Algorithm~\ref{alg:train}                          &0.9 * 121848.27 \tnote{2} \\
        num\_prelabelled            &Algorithm~\ref{alg:train} / mix\_in\_prelabelled\_data &8192\\
        ratio\_test                 &Algorithm~\ref{alg:train} / split\_data            &0.1\\
        learning\_rate              &Algorithm~\ref{alg:train} / Adam                   &0.001\\
        weight\_decay               &Algorithm~\ref{alg:train} / Adam                   &5e-06 \\
        patience                    &Algorithm~\ref{alg:train} / ReduceLROnPlateau      & 2\\
        factor                      &Algorithm~\ref{alg:train} / ReduceLROnPlateau      & 0.5\\
        
        \bottomrule
    \end{tabular}
    \begin{tablenotes}
     \item[1] 1022.5 was the test loss after initial discriminator training.
     \item[2] 121848.27 represents the average \predobjfixedapproxmip value over our artificial data.
   \end{tablenotes}
    \end{threeparttable}
    \caption{Parameters for training.}
    \label{tab:initial_discriminator}
\end{table}


%% file: sections/appendix/MIPParams.tex
\begin{table}[ht]
    \centering
    \begin{tabular}{ll}
        \toprule
        Parameter                   & Value \\
        \midrule
        TimeLimit                   &  3600 (s) \\
        FeasibilityTol              &  1e-6 \\
        MIPGap                      &  1e-4 \\
        MIPGapAbs                   &  1e-2 \\
        NumericFocus                &  3 \\
        
        \bottomrule
    \end{tabular}
    \caption{Parameters for MIP solving.}
    \label{tab:mip_params}
\end{table}

%% file: sections/appendix/parameterTable.tex
\begin{table}
    \centering
    \begin{tabular}{llll}
\toprule
{} & Parameters & Inception Blocks & Small Inception Blocks\\
\midrule
Neural Network          &               1,701,505 &                      13 &                              12 \\
Generator               &               1,165,576 &                      13 &                               0 \\
Discriminator           &                 535,929 &                       0 &                              12 \\
Inception Block         &                  87,296 &                       - &                               - \\
Small Inception Block   &                  27,936 &                       - &                               - \\
\bottomrule
\end{tabular}
    \caption{Number of parameters in the neural network and submodules.}
    \label{tab:nn_architecture}
\end{table}